\documentclass[10pt]{article}
 \pdfoutput=1
\usepackage[margin=1in]{geometry}
	 \usepackage[english]{babel}
	 \usepackage{amsmath}
	 \usepackage{amsfonts}
	 \usepackage{amssymb,multicol}
	 \usepackage{graphicx}
	 \usepackage{epstopdf}	
	 \usepackage{hyperref}
	 \usepackage{datetime}
	 \usepackage{ulem}
	  \usepackage{xspace}
	 \usepackage{bm}
	 \usepackage{lscape}
	 \usepackage{lineno,hyperref}
	 \usepackage {color}
	 \usepackage{mathtools}
	 \modulolinenumbers[5]

\definecolor{lightgray}{gray}{0.8}
\definecolor{llightgray}{gray}{0.95}

\def\##1{{\bf #1}}
\def\=#1{\underline{\underline{#1}}}

\def\*#1{\ensuremath{\breve{{#1}}}}
\def\+#1{\ensuremath{\breve{\underline{\underline #1}}}}

\renewcommand{\~}[1]{\bm{{{\bar{#1}}}}}
\renewcommand{\#}[1]{\textbf{#1}}

\def\.{\mbox{ \tiny{$^\bullet$} }}

\def\eps{\varepsilon}

\def\epsrel{\eps_{\rm rel}}
\def\epsinvrel{\beta_{\rm rel}}
\def\rel{_{\rm rel}}

\def\O{\scriptscriptstyle 0}
\def\epso{\eps_{\scriptscriptstyle 0}}

\def\lambdao{\lambda_{\scriptscriptstyle 0}}
\def\muo{\mu_{\scriptscriptstyle 0}}

\def\etao{\eta_{\scriptscriptstyle 0}}

\def\ko{k_{\scriptscriptstyle 0}}
\def\co{c_{\scriptscriptstyle 0}}
\def\Eo{E_{\, \scriptscriptstyle 0}}

\def\lambdaomin{\lambda_{\O_{min}}}
\def\lambdaomax{\lambda_{\O_{max}}}

\def\sfE{{\sf E}}
\def\sfEc{{\sf E}_{\rm c}}
\def\sfEFn{{\sf E}_{\rm F_{\rm n}}}
\def\sfEFp{{\sf E}_{\rm F_{\rm p}}}
\def\sfEg{{\sf E}_{\rm g}}
\def\sfEi{{\sf E}_{\rm i}}
\def\sfEv{{\sf E}_{\rm v}}

\def\zetag{\zeta_{\rm g}}

 \def\sfJ{{\sf J}}
\def\sfZ{{\sf Z}}
\def\sfm{{\sf m}}

\def\sff{{\sf f}}
\def\sfc{{\sf c}}

\def\sfvm{{\sf v}_{\sf m}}

\def\mfv{{\underline v}}
\def\mfM{\=M}

\def\zero{^{\circ}}

\def\ux{\hat{\underline x}}
\def\uy{\hat{\underline y}}
\def\uz{\hat{\underline z}}
  
 \def\JscOpt{\ensuremath{J_{\rm SC}^{\rm Opt}}}
\def\Jsc{\ensuremath{J_{\rm SC}}}

\def\Voc{\ensuremath{V_{\rm OC}}}

\newcommand{\ipt}[1]{<#1>_{\gamma}}
\newcommand{\pt}[1]{\left(#1\right)_{\gamma}}

\newcommand{\dsbg}[1]{\left[\!\left[ #1 \right]\!\right]_{\gamma}}

\def\ap{a_{\rm p}}
\def\as{a_{\rm s}}

\def\Edc{E_{\rm dc}}

\def\Gs{G_{\rm s}}

\def\Jn{J_{\rm n}}
\def\Jp{J_{\rm p}}
\def\Js{J_{\rm s}}

\def\kB{k_{\rm B}}

\def\Lp{L_{\rm p}}
\def\Li{L_{\rm i}}
\def\Lg{L_{\rm g}}
\def\Lm{L_{\rm m}}
\def\Ln{L_{\rm n}}
\def\Ls{L_{\rm s}}
\def\Lt{L_{\rm t}}
\def\Lx{L_{\rm x}}
\def\Lz{L_{\rm z}}

\def\Nc{N_{\rm c}}
\def\Nf{N_{\rm f}}

\def\NP{N_{\rm P}}
\def\Ns{N_{\rm s}}
\def\Nt{N_{\rm t}}
\def\Nv{N_{\rm v}}
\def\Nz{N_{\rm z}}

\def\Omegaph_{{\Omega_{\rm ph}}}
\def\Omegael_{{\Omega_{\rm el}}}

\def\Pdeg{P_{\rm deg}}
\def\Pin{P_{\rm in}}
\def\Pmax{P_{\rm max}}

\def\Ideg{I_{\rm deg}}

\def\Vext{V_{\rm ext}}
\def\Vth{V_{\rm th}}

\def\phin{\phi_{\rm n}}
\def\phip{\phi_{\rm p}}
\def\phis{\phi_{\rm s}}

\def\mun{\mu_{\rm n}}
\def\mup{\mu_{\rm p}}
\def\mus{\mu_{\rm s}}

\def\nL{n_{\rm L}}
\def\nR{n_{\rm R}}

\def\pL{p_{\rm L}}
\def\pR{p_{\rm R}}

\def\qe{q_{\rm e}}

\def\taun{\tau_{\rm n}}
\def\taup{\tau_{\rm p}}

\def\dz{d_{\rm z}}

\begin{document}

\begin{center}

\textbf{Coupled Optoelectronic Simulation and Optimization of Thin-Film Photovoltaic Solar Cells}
\vspace{4mm}

 {Tom H. Anderson}, {Benjamin J. Civiletti}, {Peter B. Monk}\\
 \textit{University of Delaware, Department of Mathematical Sciences,
	Newark, DE 19716, USA}
{Akhlesh Lakhtakia}\\
\textit{Pennsylvania State University, Department of Engineering Science and Mechanics, University Park, PA 16802, USA}
\end{center}

\begin{abstract}
{A design tool was formulated for optimizing 
the efficiency of inorganic, thin-film, photovoltaic solar cells. The solar cell
can have multiple semiconductor layers in addition to antireflection coatings, passivation layers, and buffer layers. The solar cell is backed by a metallic grating which is periodic along a fixed
direction. The rigorous coupled-wave approach  is used to calculate the electron-hole-pair
generation rate. The hybridizable discontinuous Galerkin method is 
used to solve the drift-diffusion equations that govern  charge-carrier transport in the semiconductor layers. The chief output is the solar-cell efficiency which is maximized using
 the differential evolution algorithm to determine the optimal dimensions
and bandgaps of the semiconductor layers.}\\

\end{abstract}

\vspace{3mm}

\section{Introduction}
	
The simulation of thin-film photovoltaic solar cells requires {the coupling of (i)}
an optical model capable of {capturing the absorption of photons with (ii) an electrical model capable of simulating the  transport of   charge carriers throughout the solar cell \cite{Jenny_book,Fonash_book}}. To optimize solar cell designs, both models needed to be tailored to be rapidly computable,  accurate, and  robust across the
relevant parameter space. To this end, {we have developed a coupled optoelectronic simulation
technique for thin-film photovoltaic solar cells containing periodic structures, and we have used this simulation technique  together} with the differential evolution algorithm (DEA)~\cite{Storn,Swagatam2011} to optimize device performance.

The solar cell is assumed to be infinitely extended in the $xy$-plane,
periodic in the $x$-direction, translation invariant in the
$y$-direction, and of finite extent in the $z$-direction.  These assumptions allow the simulation to be
limited to two-dimensional (2D) space.  However, we emphasize
that our approach will work in three-dimensional (3D)
space at the cost of greater computational complexity and cost.
	
The first step of the coupled optoelectronic simulation involves modeling the photonic characteristics of the solar cell. In the photonic step, the rigorous coupled-wave approach (RCWA)~\cite{Lalanne96,Lalanne97,Polo_book} is used to determine the electromagnetic
fields in the solar cell due to incident solar radiation. The RCWA is an efficient computational technique
to solve the frequency-domain Maxwell equations and determine the electric field in a representative subsection of the periodic (2D or 3D) domain.  {The  frequency-domain electric field  is used to 
determine the
photon absorption rate, and therefore the generation rate of electron--hole pairs,}
  throughout the solar cell. Rapid calculation
 across the solar spectrum \cite{Civiletti,Faiz18} is possible due to RCWA being a pseudo-spectral method, i.e., it avoids spatial meshing in the periodic dimensions.

The second step  of the coupled optoelectronic simulation involves modeling the electronic characteristics of the solar cell. In the electronic step, the
electron-hole-pair generation rate is used as the input to
a drift-diffusion electronic (DDE) model in order  to calculate  the current-voltage curve
of the solar cell and, hence, the efficiency of the solar cell.
Because of the cost of solving the DDE model, we  average the generation rate across one period in the
$xy$-plane and implement the DDE model for
variations of the electric potential and charge-carrier transport along the $z$-axis.  The resulting one-dimensional (1D) DDE model consists of six coupled differential equations.   The nonlinear Shockley--Read--Hall, Auger, and radiative  terms are included  to model electron-hole recombination \cite{Jenny_book,Fonash_book}. A hybridizable discontinuous Galerkin (HDG) method
\cite{Brinkman,Lehrenfeld,CockburnHDG,FuQiuHDG} is developed, and the Newton--Raphson method
\cite{Jaluria_book} is used to find a solution of the resulting nonlinear system.   

A version of the HDG method has already been used to 
simulate organic solar cells~\cite{Brinkman},   based
on the work of Lehrenfeld~\cite{Lehrenfeld}.  We use 
a similar scheme based on the classical HDG method~\cite{CockburnHDG}, following Fu 
\textit{et al.}~\cite{FuQiuHDG} who have analyzed the method for a linear convection-diffusion system.  Our study points to the advantages of this method in {simulating \textit{inorganic} solar cells having a spatially variable electron-hole-pair generation rate in the presence of semiconductor}
heterojunctions.

Anderson \textit{et al.} \cite{Anderson17} recently reported  a coupled optoelectronic simulation
of a Schottky-barrier thin-film solar cell, with the finite-element method (FEM) \cite{ChenFEM_book} used for both the photonic and
the electronic steps of the simulation. Later, Anderson \textit{et al.} \cite{Anderson18} used the RCWA for
the photonic step and the standard FEM for the electronic step for simulating 
 the same Schottky-barrier thin-film solar cell (which is not considered in this paper). Other researchers have also developed techniques for
 optoelectronic simulation \cite{KrcTopic,Saiprasad,solcore}.
 
Our emphasis in this paper is different in two ways: first, our simulations lead
to optimization of  thin-film photovoltaic solar  cells containing periodic structures;
second, we
concentrate on understanding the numerical performance of each of the two steps
in the coupled optoelectronic simulation.  For example, we investigate the convergence rate of RCWA for two model problems.  We also propose, test, and use
a high-order HDG scheme using fifth-degree polynomials for the DDE model, 
quite possibly the first time that such a high-order scheme has been used.  Since the use of DEA
for optimization of solar cells has been investigated elsewhere by us
\cite{Civiletti,Anderson18}, we devote little space to it in this paper.

The photonic and electronic steps described in this paper are closely related to components of 
the  Solcore software package described  recently by
Alonso-\'{A}lvarez \textit{et al.} \cite{solcore}.  That software package includes an RCWA solver
as well as a 1D DD solver based on quasi-Fermi levels.  
{In contrast, our DD solver} uses the density of holes and electrons and is based on a 
high-order HDG scheme.

The ultimate goal is to maximize the efficiency of the solar cell, and this goal distinguishes our work from that of Alonso-\'{A}lvarez \textit{et al.}~\cite{solcore}.  Because of the presence of local minima, we use the DEA~\cite{Storn,Swagatam2011}
for optimization.  This requires the evaluation of the
efficiency of solar-cell designs from  all regions of the multi-dimensional space encompassing appropriate parameters,
which thus necessitates the development of a robust and efficient solver.   See
\href{https://www.pvlighthouse.com.au}{https://www.pvlighthouse.com.au}
for a list of other available photovoltaic-design software. Our work here should be seen as a complement to those efforts by extending
the use of coupled optoelectronic simulation to the optimal design of solar cells as well as testing new simulation techniques.

\begin{figure}
\begin{center}
\begin{tabular}{cc}
\resizebox{0.49\textwidth}{!}{\includegraphics{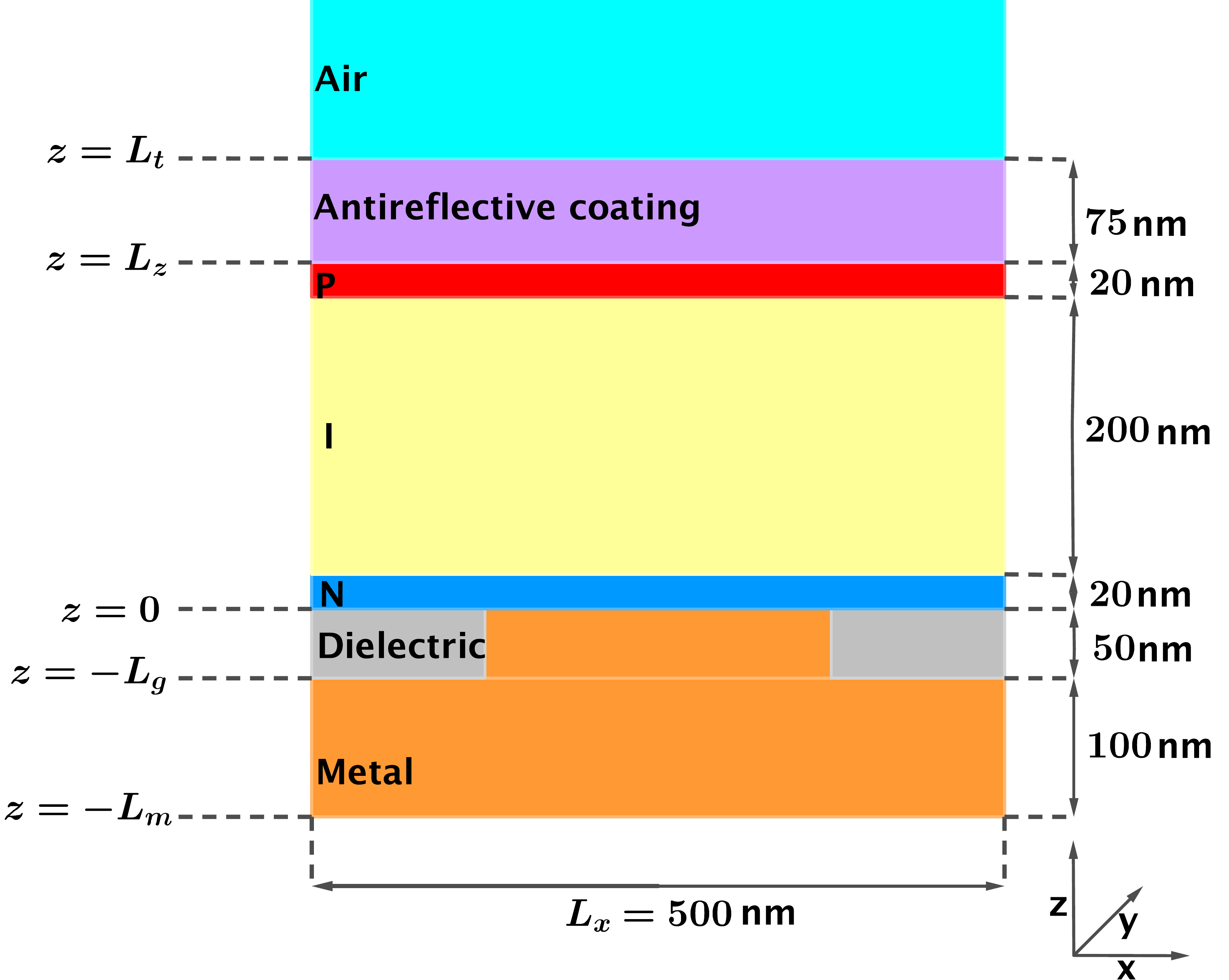}}&
\resizebox{0.48\textwidth}{!}{\includegraphics{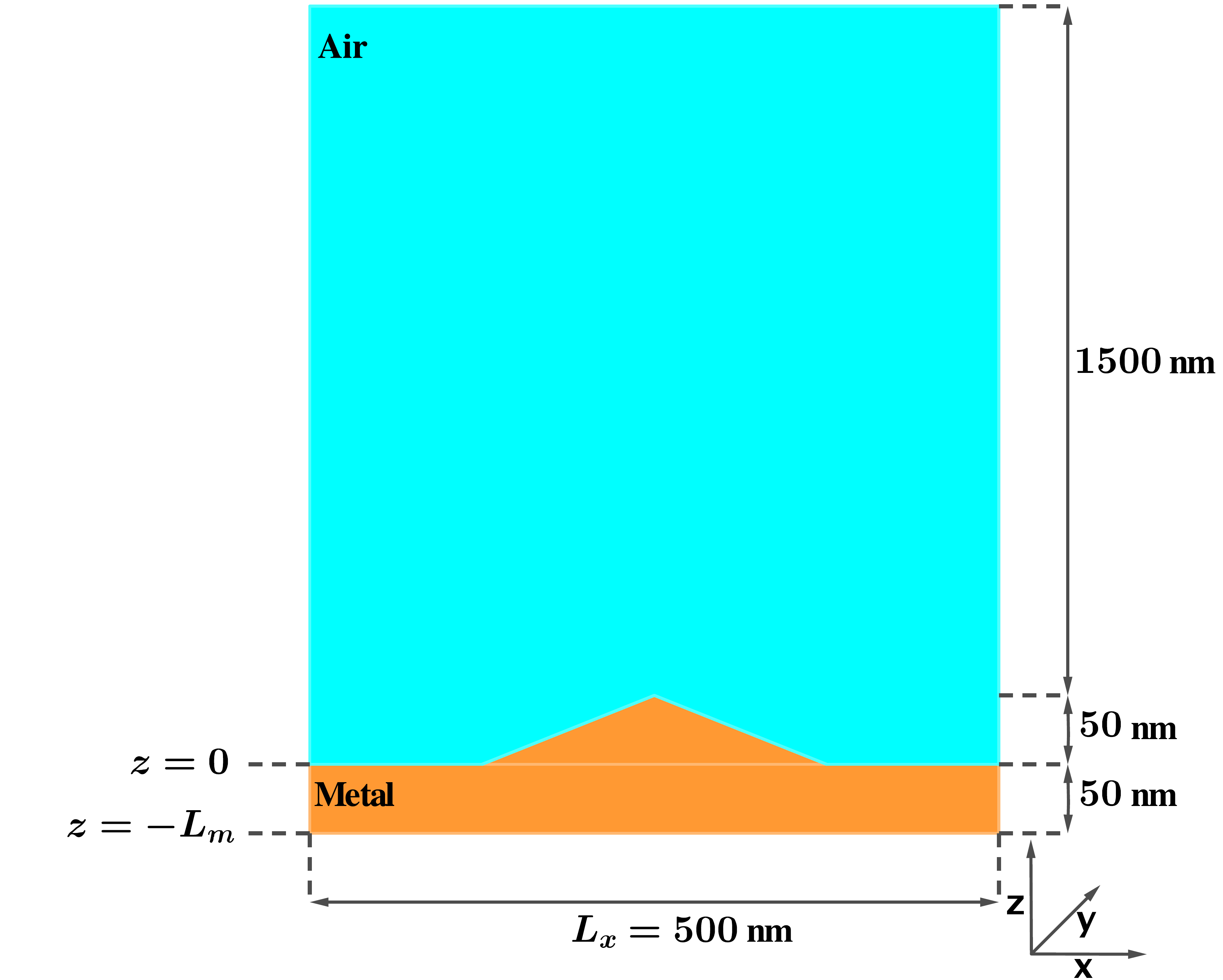}}\\(a)&(b)
\end{tabular}
\end{center}
\caption{(a) An academic model of a metal-backed \textit{p-i-n} solar cell 
used as an example throughout this paper.  The dimensions shown are simply for illustration. (b) A  triangular metal grating in air used to test RCWA (of course, this is not a solar cell). Figures are not to scale.}
\label{Ben1}
\end{figure}

In this paper, we assume that the solar cell occupies the region 
\begin{equation}
{\cal X}:
\{\underline{r}\equiv(x,y,z)\;|\; -\infty<x<\infty,\;-\infty<y<\infty,\; -\Lm<z<\Lt\}
\end{equation}
in $\mathbb{R}^3$, with
air occupying the region $\mathbb{R}^3\setminus{\cal X}$. As stated previously,
the device is assumed to be periodic along the $x$-direction with period $\Lx$  and translation invariant 
along the  $y$-direction. Solely for illustrating various simulation issues, we performed calculations for
the device shown Fig.~\ref{Ben1}(a). The current-generating region of this device comprises
a thin layer of $p$-type semiconductor, a thicker layer of $i$-type semiconductor,
and a  thin layer of $n$-type semiconductor \cite{Jenny_book,Fonash_book}. 
The plane $z=0$ is taken to be the bottom face of the $n$-type semiconductor
and  the plane $z=\Lz$ is taken to be the top face of the $p$-type semiconductor.
 Below  the plane $z=0$ is a grating formed by
infinitely long strips of a metal and a dielectric material, the grating period being $\Lx$. Below the grating
is a metal layer that is sufficiently thick to prevent light from escaping into air below it. Above the
$p$-type semiconductor is an antireflection coating.  The plane $z=\Lt$ is the top face
of the antireflection coating and the plane $z=-\Lm$ is the bottom face of the metal layer below the grating.
The grating shape in Fig.~\ref{Ben1}(a) is an academic example. A triangular grating used to test RCWA  is shown in Fig.~\ref{Ben1}(b).

The layout of the paper is as follows.  In Section~\ref{Omodel}, we start by describing the optical model
solved in the photonic step, summarize the 2D RCWA, and define the 
charge-carrier generation rate that is used as the input to the DDE model solved
in the electronic step.  In Section~\ref{Sec:RCWAConv} we show two examples of RCWA
convergence to illustrate the advantages and drawbacks of this approach for solar
cells. Section~\ref{Emodel} is devoted to the electronic step. The DDE model
is described in Section~\ref{DD}. The boundary conditions  derived from local quasi-thermal equilibrium (LQTE)  are discussed
in Section~\ref{LQTE}, while the electron-hole recombination terms are presented in Section~\ref{Recom}. We discuss
heterojunctions in Section~\ref{Sec:Hetero}, and summarize a non-dimensionalized DDE model in Section~\ref{dimless}.
Then we proceed to describe the HDG method in Section~\ref{HDG}, starting with 
the formulation of a generalized transport system in Section~\ref{sec:gts}, followed by
the discretization
of 1D space (along the $z$ axis)  in Section~\ref{HDGD} and
implementational aspects of the convection term in Section~\ref{Sec:Nonlinear}.  Next, in Section~\ref{upwind}, we present the upwinding strategy needed to stabilize the DDE model, followed in Section~\ref{numrecom} by the interpolation method for implementing the non-linear recombination terms. {The implementation of jump conditions is presented
in Section~\ref{numhetero}.} The numerical model is completed by discussing the homotopy (or continuation) method for solving the nonlinear system in Section~\ref{homotopy}.  A numerical test is presented in Section~\ref{HDGtest}. Together, the photonic and electronic steps provide information used by the DEA optimization scheme, which is described briefly in Section~\ref{DEA}. The paper culminates in a presentation of results for the model solar cell presented in Fig.~\ref{Ben1}(a).  We end with some conclusions in Section~\ref{concl}.

A note about notation: Underlined quantities represent vectors, with $\hat{\underline{x}}$, $\hat{\underline{y}}$ and {$\hat{\underline{z}}$} being the unit vectors along the $x$-, $y$-, and ${z}$-axes. The $\;\*{}\;$ mark is used to denote quantities emerging from the implementation of the  spatial Fourier transform
\cite{Goodman} with respect to
$x$. For the frequency-domain calculations carried out in the photonic step
(Section~\ref{Omodel}), electromagnetic fields
are taken to vary with time $t$ harmonically as $\exp(-i\omega t)$, where $\omega$ is the angular frequency and $i=\sqrt{-1}$. Air is assumed to have the same electromagnetic properties as free space or vacuum. The wavenumber in air is denoted by
$\ko=\omega/\co$ where $\co$ is the speed of light in air, and $\lambdao=2\pi/\ko$
is the wavelength of light in air.  The {permittivity and permeability of free space are denoted, respectively, by $\epso=8.854188\times10^{-12}$~F~m$^{-1}$ and $\muo=4\pi\times10^{-7}$~H~m$^{-1}$, and $\etao=\sqrt{\muo/\epso}$} is the intrinsic impedance of free space.

\section{Photonic Step} \label{Omodel}
In this section, first we   summarize the RCWA for solving the frequency-domain
Maxwell equations, and then we present some numerical results to quantify the performance of  RCWA,
which is the workhorse
for computational investigations of optical gratings \cite{Maystre,Moharam}. The adequacy
of this approach for solar cells has been established by comparison to other approaches
\cite{SolanoFEM,shuba15,lokar18}. Oddly, the numerical convergence of this method has not been investigated before; hence, a part of this section is devoted to a numerical investigation of RCWA
convergence.

\subsection{RCWA Formulation}
The solar cell is taken to be {illuminated from the half space $z>\Lt$
by a normally incident plane wave whose electric field phasor is
 denoted by
	\begin{equation}
	{\underline E}_{\rm inc}(x,z,\lambdao)=\Eo\frac{\ap\, \ux+\as\,\uy}{\sqrt{\ap^{2}+\as^{2}}}\exp(-i\ko z).
	\end{equation}
The parameters $\as$ and $\ap$ determine the polarization state of the incident light; 
thus,  $\ap = 1$ and $\as = 0$ for $p$-polarized light, but $\as=1$ and $\ap=0$ for $s$-polarized light. For solar-cell calculations, we set $\ap=\as=1$ because direct sunlight can be assumed to be unpolarized.}

Due to the periodicity of the solar cell in the $x$-direction, the  optical relative permittivity $\epsrel(x,z,\lambdao)$  is represented everywhere by the Fourier series
\begin{align}
\epsrel(x,z,\lambdao) &= \sum_{\ell=-\infty}^{\infty} {\eps}\rel^{(\ell)}(z,\lambdao)\exp\left(i\kappa^{(\ell)}x\right),\quad \vert z\vert <\infty,\quad  \vert x\vert <\infty,
\label{eq3}
\end{align}
where ${\eps}\rel^{(\ell)}(z,\lambdao)$ are   Fourier coefficients and  $\kappa^{(\ell)} =   2\pi \ell/\Lx$.In a similar way, the optical relative impermittivity $\epsinvrel(x,z,\lambdao)=1/\epsrel(x,z,\lambdao)$
is repesented everywhere by the Fourier series
 \begin{align}
\epsinvrel(x,z,\lambdao) &= \sum_{\ell=-\infty}^{\infty} {\beta}\rel^{(\ell)}(z,\lambdao)\exp\left(i\kappa^{(\ell)}x\right),\quad \vert z\vert <\infty,\quad  \vert x\vert <\infty,
\label{eq4}
\end{align}
where ${\beta}\rel^{(\ell)}(z,\lambdao)$ are   Fourier coefficients.	We can also express the $x$-dependences of the electric and magnetic field phasors by their Fourier series as 
\begin{align}
{\underline E}(x,z,\lambdao) &= \sum_{\ell=-\infty}^{ \infty} {\underline e}^{(\ell)}(z,\lambdao)\exp\left(i\kappa^{(\ell)}x\right),\quad \vert z\vert <\infty,\quad  \vert x\vert <\infty,\\
{\underline H}(x,z,\lambdao) &= \sum_{\ell=-\infty}^{ \infty} {\underline h}^{(\ell)}(z,\lambdao)\exp\left(i\kappa^{(\ell)}x\right), \quad \vert z\vert <\infty,\quad  \vert x\vert <\infty,
\end{align}
where ${\underline e}^{(\ell)} \equiv {e}_x^{(\ell)}\ux+{e}_y^{(\ell)}\uy+{e}_z^{(\ell)}\uz$ and ${\underline h}^{(\ell)}\equiv {h}_x^{(\ell)}\ux+{h}_y^{(\ell)}\uy+{h}_z^{(\ell)}\uz$ are   Fourier coefficients and normal incidence is implicit.

For computational tractability, all four of the foregoing series are truncated to include only  $\ell \in \left\{-\Nt,...,\Nt\right\}$, $\Nt\geq0$. The Fourier coefficients are collected in the $(2\Nt+1)$-column vectors
\begin{eqnarray}
\nonumber
&&{\*\eps}\rel(z,\lambdao)=\left[ 
\eps\rel^{(-\Nt)}(z,\lambdao),\eps\rel^{(-\Nt+1)}(z,\lambdao),...,\eps\rel^{(\Nt-1)}(z,\lambdao),\eps\rel^{(\Nt)}(z,\lambdao)\right]^T,
\\
&&
\\
\nonumber
&&
{\*\beta}\rel (z,\lambdao)=\left[ 
\beta\rel^{(-\Nt)}(z,\lambdao),\beta\rel^{(-\Nt+1)}(z,\lambdao),...,\beta\rel^{(\Nt-1)}(z,\lambdao),\beta\rel^{(\Nt)}(z,\lambdao)\right]^T,
\\
&&
\\
\nonumber
&&
{\*e}_\sigma(z,\lambdao)=\left[ 
e_\sigma^{(-\Nt)}(z,\lambdao),e_\sigma^{(-\Nt+1)}(z,\lambdao),...,e_\sigma^{(\Nt-1)}(z,\lambdao),e_\sigma^{(\Nt)}(z,\lambdao)\right]^T,
\\
&& \hspace{6cm}\sigma\in\{x,y,z\}\,,
\\
&&
\nonumber
\*h_\sigma(z,\lambdao)=\left[ h_\sigma^{(-\Nt)}(z,\lambdao),h_\sigma^{(-\Nt+1)}(z,\lambdao), ...,h_\sigma^{(\Nt-1)}(z,\lambdao),h_\sigma^{(\Nt)}(z,\lambdao)\right]^T
\\
&& \hspace{6cm}\sigma\in\{x,y,z\}\,,
\end{eqnarray}
where the superscript $T$ denotes the transpose.   Furthermore,   
  the  $(2\Nt+1)\times (2\Nt+1)$ matrix  
\begin{align}
\*K= 
{\rm diag}\left[  \kappa^{(-\Nt)},\, \kappa^{(-\Nt+1)},\,..., 
\kappa^{(\Nt-1)},\,\kappa^{(\Nt)}\right]
\end{align}  
is defined for convenience. Finally, the
Toeplitz matrix of the Fourier coefficients $\{\xi^{(\ell)}\}_{\ell=-\infty}^{\infty}$ of a periodic function $\xi(x)=\xi(x\pm \Lx)$ is defined as
\begin{align}
{\cal{}T}_{\Nt}(\xi)=
\left[
\begin{array}{ccccc}
\xi^{(0)} & \xi^{(-1)} &  ... &  \xi^{(-2\Nt+1)} & \xi^{(-2\Nt)}\\[5pt]
\xi^{(1)} & \xi^{(0)} &   ... &   
\xi^{(-2\Nt+2)} &\xi^{(-2\Nt+1)}\\[5pt]
\vdots & \vdots & \ddots  &\vdots & \vdots \\[5pt]
\xi^{(2\Nt-1)} & \xi^{(2\Nt-2,\lambdao)} &   ... &   
\xi^{(0)} & \xi^{(-1)}\\[5pt]
\xi^{(2\Nt)} & \xi^{(2\Nt-1)} &  ... &   
\xi^{(1)} & \xi^{(0)}
\end{array}
\right]\,. \label{Eq2.105}
\end{align}

Substitution of the four Fourier series into the frequency-domain Maxwell curl equations yields two
matrix ordinary differential equations \cite[Sec.~2.3.4]{Polo_book}.
The equation
\begin{eqnarray}
\nonumber
&&
\frac{d}{dz}
\left[ \begin{array}{c}
\*e_x(z,\lambdao)\\[3pt]
\*h_y(z,\lambdao)
\end{array}
\right]
= 
\\[10pt]
&& 
 i 
\left[
\begin{array}{cccc}
\*0&\omega\muo\*I-  \left(\omega\epso\right)^{-1}\*K\*N (z,\lambdao)\*K\\[3pt]
\omega\epso\*M (z,\lambdao) &\*0
\end{array}\right]
\left[ \begin{array}{c}
\*e_x(z,\lambdao)\\[3pt]
\*h_y(z,\lambdao)
\end{array}
\right]
\label{Eqn:diff-p}
\end{eqnarray}
must be solved if $\ap\ne0$, and
the equation
\begin{eqnarray}
\nonumber
&&
\frac{d}{dz}
\left[ \begin{array}{c}
\*e_y(z,\lambdao)\\[3pt]
\*h_x(z,\lambdao)
\end{array}
\right]
= 
\\[10pt]
&& 
- i 
\left[
\begin{array}{cccc}
\*0&\omega\muo\*I \\[3pt]
\omega\epso\*M (z,\lambdao)-(\omega\muo)^{-1}\*K^2
 &\*0
\end{array}\right]
\left[ \begin{array}{c}
\*e_y(z,\lambdao)\\[3pt]
\*h_x(z,\lambdao)
\end{array}
\right]
\label{Eqn:diff-s}
\end{eqnarray}
if $\as\ne0$. 
In these equations, $\*0$ is the $(2\Nt + 1)\times(2\Nt+1)$ null matrix, and
$\*I$ is the $(2\Nt + 1)\times(2\Nt+1)$ identity matrix. 
As a corollary to the mutual independence of  Eqs.~(\ref{Eqn:diff-p}) and (\ref{Eqn:diff-s}),
 $p$-polarized light does not interact with $s$-polarized light and \textit{vice versa}.

We solved both 
Eqs.~(\ref{Eqn:diff-p}) and (\ref{Eqn:diff-s}) because direct sunlight can be assumed to be unpolarized.
Thereafter,
 the Fourier coefficients of the $z$-components of the electric
 and magnetic field phasors were obtained as
\begin{align}
\left.\begin{array}{l}
\*{e}_z(z,\lambdao) = -( \omega\epso)^{-1} \*N (z,\lambdao) 
\#{\*{K}}  \*{h}_y(z,\lambdao)\\[5pt]
\*{h}_z(z,\lambdao) =(\omega \muo)^{-1}  \#{\*{K}} \*{e}_y(z,\lambdao)
\end{array}\right\}.\label{Eqn:ezhz}
\end{align}

The matrices $\*M $ and $\*N $ 
in Eqs.~(\ref{Eqn:diff-p}) and (\ref{Eqn:diff-s}) 
encode the Fourier series provided in Eqs.~(\ref{eq3})
and (\ref{eq4}) for $\epsrel$ and $\epsinvrel$.
There are four choices possible for the combination
of $\*M $ and $\*N $,  the best
choice depending on the polarization state of the incident plane wave~\cite{li96}. 
Following  Lalanne and Morris~\cite{Lalanne96}, we first define the Toeplitz matrices
\begin{equation}
\label{eq16}
\*{E}(z,\lambdao)= \mathcal{T}_{\Nt}[\*{\eps}\rel(z,\lambdao)]
\end{equation}
and
\begin{equation}
\label{eq17}
\*{B}(z,\lambdao)= \mathcal{T}_{\Nt}[\*{\beta}\rel (z,\lambdao)].
\end{equation}
Then the four possible choices are identified as follows:
\begin{description}
\item[Choice 1:]  The most \textit{direct} choice is for $\*M $ to encode $\epsrel$ and $\*N $ to encode
$\eps\rel^{-1}$ so that
$\*M   =\*E$
and
$\*N   =\*B$.  
\item[Choice 2:]  We can write $\epsrel =1/\epsinvrel$ and arrive at the choice
$\*M   =\*B^{-1}$
and
$\*N   =\*B$.
\item[Choice 3:]  In contrast to Choice 2, we can write $\epsinvrel=1/\eps\rel$
to choose
$\*M   =\*E$
and
$\*N   =\*E^{-1}$.
\item[Choice 4:] The final choice  is
$\*M   =\*B^{-1}$
and
$\*N   =\*E^{-1}$, in direct contrast to Choice 1.
\end{description}
Choices 1 and 3 are identical on the one hand, and Choice 2 is the same as Choice 4 on the
other, for   $s$-polarized illumination. All four choices are distinct for $p$-polarized illumination.

In order to solve Eqs.~(\ref{Eqn:diff-p}) and (\ref{Eqn:diff-s}), it is usual to discretize in the $z$-direction by subdividing  the solar cell into thin   slices of infinite extent in the $xy$-plane~\cite{Polo_book}.  In each slice,
$\epsrel$ is uniform in the $z$-direction, but it is either uniform or piecewise uniform in the $x$-direction. Depending on the profile of the grating, this may require  piecewise-uniform
approximations of $\epsrel$ and $\epsinvrel$. 
Equations~(\ref{Eqn:diff-p}) and (\ref{Eqn:diff-s}) can then be solved exactly in each slice
 using a stabilized stepping algorithm ~\cite{Polo_book,Moharam} akin to Gaussian elimination without pivoting \cite{Jaluria_book}, so that 
 $\*{e}_x$, $\*{e}_y$ , $\*{h}_x$,  and $\*{h}_y$  are determined everywhere in the solar cell.
Equations~(\ref{Eqn:ezhz}) then let us determine $\*{e}_z$ and $\*{h}_z$  throughout the solar cell.

The convergence rates for the Fourier series of $\epsrel$ and $\epsinvrel$ are limited by the discontinuities in the relative permittivity in the slices in the grating region.  Additionally, approximating the interfaces of general grating profiles using stairstepping can also give low-order convergence with respect to the slice
thickness.  However, if the grating is in fact a superposition of a finite number of square waves,
the only approximation is due to the truncation of the Fourier series. 

\subsubsection{Electron-Hole-Pair Generation Rate}
With the assumption that the absorption of every photon in any semiconductor layer excites
an electron-hole pair, the electron-hole-pair generation rate in the semiconductor region
{$\left\{\vert{x}\vert<\Lx/2, 0 < z < \Lz\right\}$} can be calculated as
	\begin{equation}
	\label{def-Gxz}
	G(x,z)= \frac{1}{\hbar\co }\int_{\lambdaomin}^{\lambdaomax}
	\text{Im}\left\{\epsrel(x,z,\lambdao) \right\} \left\vert
	\frac{{\underline E}(x,z,\lambdao)}{E_o}
	\right\vert^2
	S(\lambdao) \operatorname{d}\!\lambdao\,,
	\end{equation}
where $\hbar$ is the reduced Planck constant and $S(\lambdao)$ is the incident power spectrum. We chose $S(\lambdao)$ to be the standard AM1.5G spectrum \cite{AM1.5}. As $S(\lambdao)$ falls off rapidly at shorter free-space wavelengths, ${\lambdaomin} = 350$~nm was fixed. Once the free-space wavelength becomes too large, the energy in a photon is not capable of exciting an electron sufficiently  to cross the bandgap. Consequently, the upper cut off was taken to be $\lambdaomax = 1240~\text{nm V} /\sfE_{\rm g,min}$, where $\sfE_{\rm g,min}$ (in V) is the minimum value of the bandgap $\sfEg$ in the solar cell.  These choices are typical values for studies on solar cells.

After assuming full quantum efficiency, i.e. each excited electron or hole is collected at the appropriate terminal, the   optical short-circuit current density can be calculated as
	\begin{align}
	\JscOpt = \qe \frac{1}{\Lx} \int_0^{\Lz}\!\!\int_{-\Lx/2}^{\Lx/2} G(x,z) \operatorname{d}\!x\operatorname{d}\!z,\label{Eqn:JscOp}
	\end{align}
where $\qe=1.60218\times10^{-19}$~C is the elementary charge. Let us note that $\JscOpt$ 
provides only a rough benchmark for the device efficiency, although it is  often used to
assess the efficiency of solar cells \cite{Dewan2009,SolanoOPT,Anttu2017}. 
However, as recombination is neglected, $\JscOpt$  is necessarily larger than the attainable short-circuit current density $\Jsc$, which is the density of current  that flows when the solar cell is illuminated and no external bias is applied (i.e., when $\Vext= 0$). We have already shown that,
when the bandgap is one of the simulation parameters, optimizing Eq.~(\ref{Eqn:JscOp}) can lead to designs with high values of $\JscOpt$ but low values of extractable electrical power \cite{Civiletti}.
 A better approach is to complement the photonic step by the electronic step, as accomplished in
Section~\ref{Emodel}.

\subsection{RCWA Convergence}\label{Sec:RCWAConv}

There are two sources of approximation error in the RCWA. The first is the truncation of the Fourier series to include only $\ell \in\{-\Nt, ..., \Nt\}$, while the second is through the discretization of space in the $z$-direction whenever there is a need to approximate the relative
permittivity  slice by slice. As it is not possible to solve Eqs.~(\ref{Eqn:diff-p}) and (\ref{Eqn:diff-s}) analytically,
in order to investigate the convergence of RCWA, the electric field phasor
${\underline E}(x,z,\lambdao)$ was calculated using an FEM code
\cite{SolanoFEM,shuba15} and compared with the  RCWA results for various values of $\Nt$, the slice thickness being fixed at a sufficiently low  value.

The FEM data were obtained using an adaptive method implemented in NGSolve~\cite{netgen}. The adaptive algorithm uses mesh bisection and the Zienkiewicz--Zhu \textit{a-posteriori} error estimator \cite{Zienkiewicz}. The FEM calculations were carried out using 5th-order continuous finite elements. Mesh adaptivity was terminated whenever the algorithm reached 100,000 degrees of freedom.
The simulated domain was sandwiched between two perfectly matched layers (PMLs). Both of the PMLs were taken to be one-wavelength thick and with a constant PML parameter of $(1.5+2.5i)$~\cite{chen03},  the reflection coefficient of the PMLs being then $\sim$$3\times 10^{-12}$.  

\subsubsection{Global Convergence of  Electric Field Phasor}
In order to investigate the global convergence of the RCWA, a simple but adequate
problem was chosen: reflection by the  grating shown in {Fig.~\ref{Ben1}(b)}. The region
above the grating was taken to be occupied by air with relative permittivity $\epsrel
= 1+10^{-9}i$, the minuscule imaginary part needed for the stability of the RCWA algorithm.
The grating, made from a fictitious metallic material with relative permittivity $\epsrel = -22+0.4i$, comprises triangular protrusions on a $50$-nm-thick film. 
Each protrusion is of height 50~nm and base 250~nm. Calculations were made with
$\lambdao=600$~nm and $\Lx= 500$~nm.

The triangular   profile was chosen for the grating, because it is approximated by a stairstep when the piecewise-uniform approximation is used to solve Eqs.~(\ref{Eqn:diff-p}) and (\ref{Eqn:diff-s}).  In addition, the true solution has a singularity in the field at the tip of the triangle which tests the algorithm's ability to handle such behavior.

As the only non-zero Cartesian component
of $\underline{E}(x,z,\lambdao)$ for $s$-polarized illumination   is the $y$-component,
 a relative global error was defined  as 
\begin{align}
e_{E_s}(\Nt) = \frac{\lVert E_{y, RCWA}(\Nt) - E_{y,FE}\rVert_2}{\lVert E_{y,FEM} \rVert_2}\label{errEs}
\end{align}
as a function of $\Nt$, with
 $E_{y, RCWA}$ denoting the  value of $E_y$ yielded by RCWA, and $E_{y,FEM}$  the value yielded
 by FEM. In addition, for any function $f(x,y)$,
\begin{align}
\lVert f \rVert_2^2 = \int_{\Omega_{\rm ph}} \vert f(x,y)\vert^2\,\mathrm{d}z\,\mathrm{d}x
\label{def-fnorm}
\end{align}
and where ${\Omega_{\rm ph}}$ is the entire simulation domain in the photonic step except the two PMLs.
For $p$-polarized illumination, the relative global error was analogously defined as
\begin{align}
e_{E_p}(\Nt) = \frac{\sqrt{\lVert E_{x, RCWA}(\Nt) - E_{x,FEM}\rVert_2^2 + \lVert E_{z, RCWA}(\Nt) - E_{z,FEM}\rVert_2^2 }}{\sqrt{\lVert E_{x,FEM} \rVert_2^2+\lVert E_{z,FEM} \rVert_2^2}},\label{errEp}
\end{align}
because only the $y$-component of $\underline{E}(x,z,\lambdao)$ is null valued.

Four different choices of $\left\{\*M,\*N\right\}$, as described after Eq.~(\ref{eq17}),
are possible when implementing the RCWA. Accordingly, four different studies
were conducted to investigate  the relative global error in the electric field phasor for 
 $s$- and  $p$-polarized illuminations.

\begin{figure}[ht]
\begin{center}
	\begin{tabular}{cc}
		(a)\includegraphics[width=0.45\textwidth]{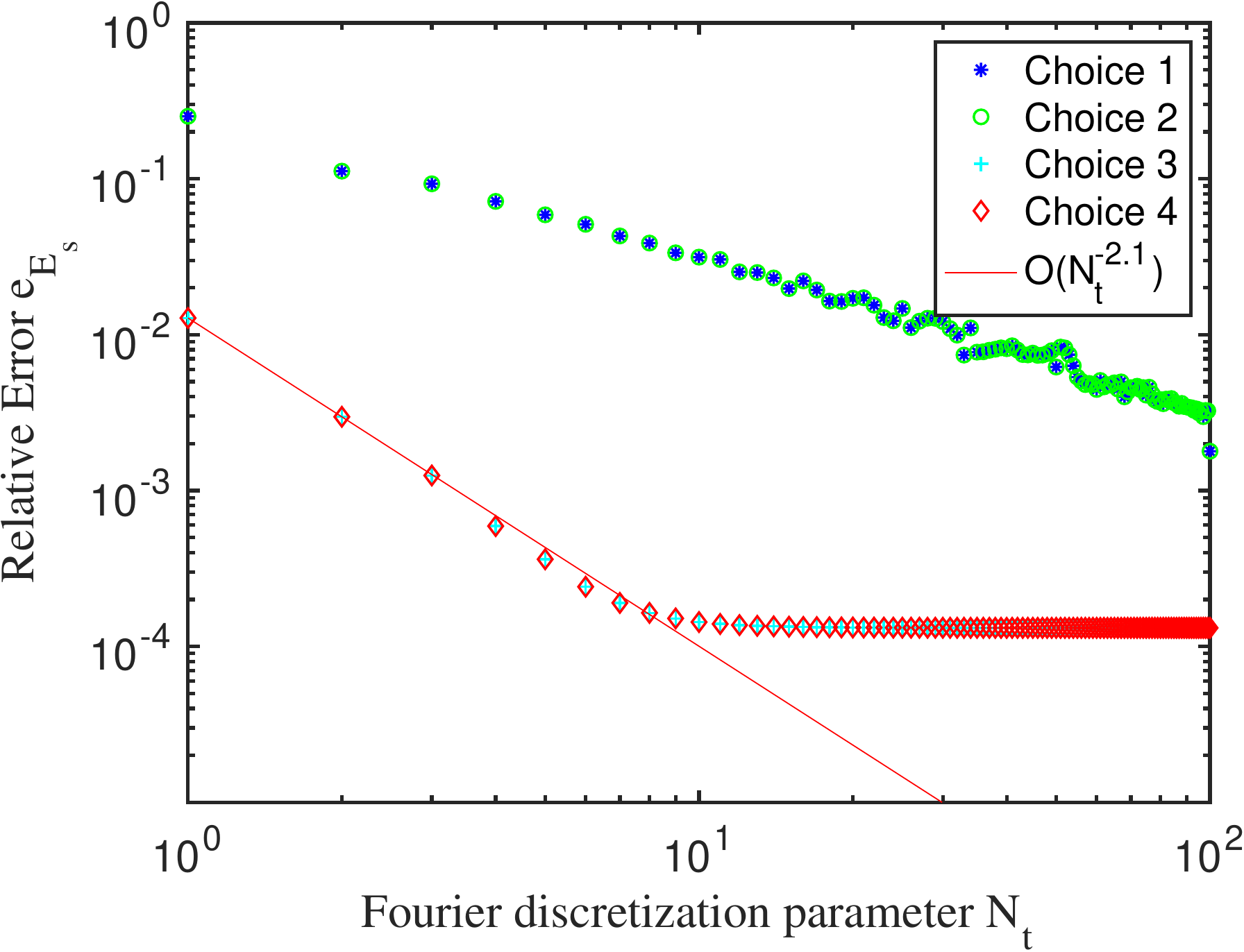}&
		(b)\includegraphics[width=0.45\textwidth]{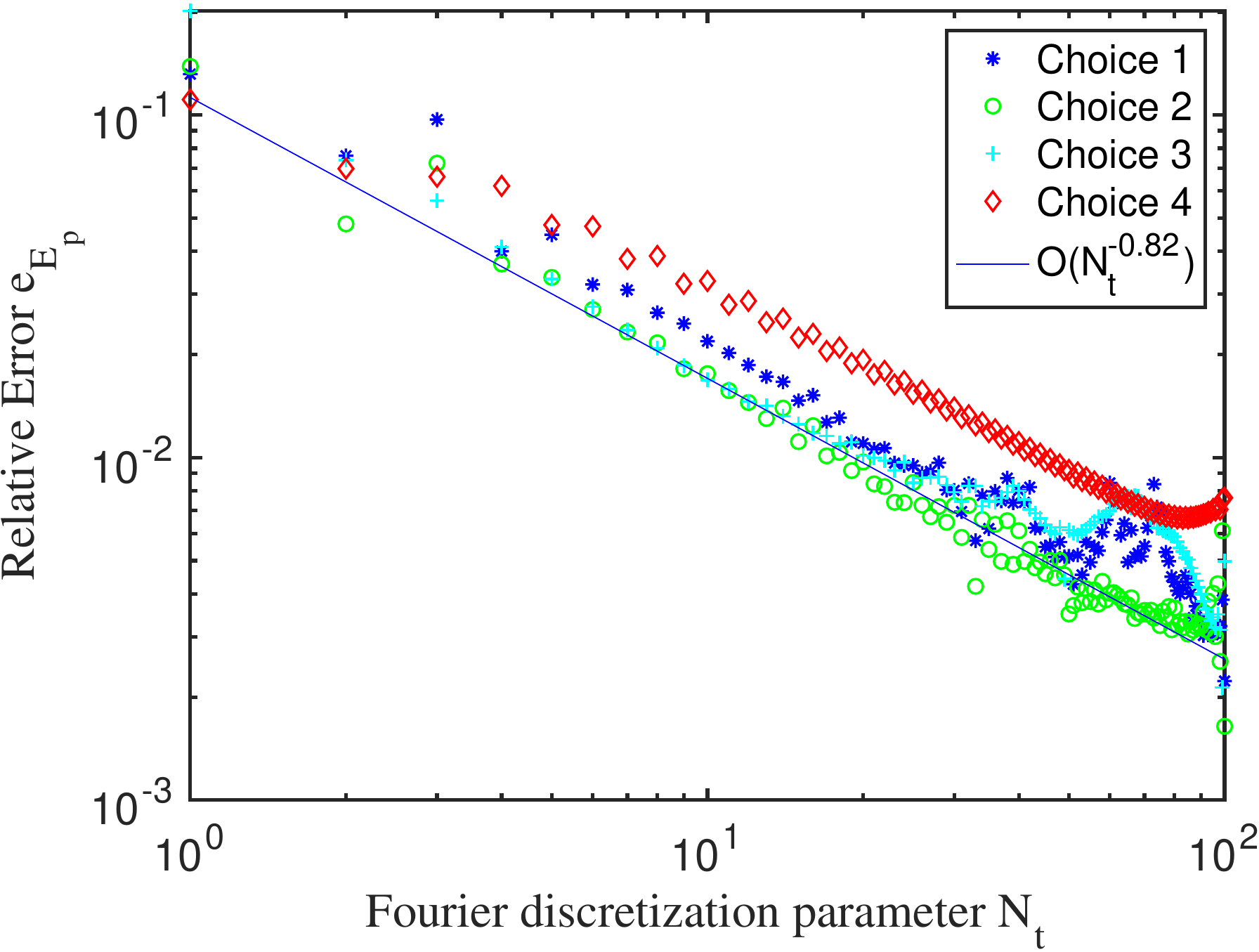}
	\end{tabular}
    \end{center}
\caption{Global relative errors (a) $e_{E_s}$ and (b) $e_{E_p}$ vs. $\Nt$ for the simple problem
depicted in {Fig.~\ref{Ben1}(b)}.
Choice 1: $\*M = \*{B}^{-1}$ and $\*N = \*{E}^{-1}$; 
Choice 2:  $\*M= \*{B}^{-1}$ and $\*N = \*{B}$; 
Choice 3: $\*M = \*{E}$ and $\*N = \*{E}^{-1}$; 
Choice 4:  $\*M = \*{E}$ and $\*N= \*{B}$.  A straight line shows the empirical order of convergence. 
}\label{Fig:RCWAConv}
\end{figure}

Plots of $e_{E_s}$ and $e_{E_p}$ vs. $\Nt$ are
presented in   Fig.~\ref{Fig:RCWAConv} for each of the four
choices of $\left\{\*M,\*N\right\}$.
For $s$-polarized illumination,  the relative global error  reduces as
$\sim{\Nt^{-2.1}}$  until {$\Nt \approx 8$} and saturates thereafter in   Fig.~\ref{Fig:RCWAConv}(a),
when $\*M = \*{E}$ (i.e., for Choices 1 and 3).  The error saturation may be due to the limited accuracy of the FEM solution {and/or} due to stairstepping error. In contrast, $e_{E_s}$ reduces as
$\sim{\Nt^{-1}}$ when $\*M= \*{B}^{-1}$ (i.e., for Choices 2 and 4).
For $p$-polarized illumination,  the relative global error decays slower than the reciprocal of $\Nt$
 in     Fig.~\ref{Fig:RCWAConv}(b),
when $\*M = \*{E}$. The decay is somewhat faster when  $\*M = \*{B}^{-1}$.

These results suggest that for calculating the electric field phasor
 due to $s$-polarized illumination, it is advantageous to choose
$\*M  = \*{E}$ because that should give second-order convergence to the true solution. 
For $p$-polarized illumination, it is preferable to choose $\*M = \*{B}^{-1}$ and $\*N = \*{E}^{-1}$, as that should give almost first-order convergence to the true solution even with stairstepping.  
The   choice of Fourier representations for
 $\left\{\*M,\*N\right\}$ thus has a profound effect on the accuracy of  the RCWA \cite{Lalanne96,li96}. 
 
\subsubsection{Convergence of  Electric Field Phasor in the Semiconductor Region}	
Correct calculation of the electric field phasor in the semiconductor region is
paramount for photovoltaic solar cells, as is clear from Eq.~(\ref{def-Gxz}).
This motivated the second convergence study using the device in {Fig.~\ref{Ben1}(a)}. 

The chosen device  contains a grating made from a  representative metal with relative permittivity $\epsrel = -22+0.4i$ and has period $\Lx= 500$~nm. The grating comprises rectangular protrusions on a $100$-nm-thick film. 
Each protrusion is of height 50~nm and base 250~nm. 
The grooves of the grating are entirely filled with a dielectric material with relative permittivity $\epsrel = 3.33 + 0.016i$.  Above this structure lies a $240$-nm-thick semiconductor region with relative permittivity $\epsrel = 9.5 + 1.25i$,  which is a representative value for a semiconductor at $\lambdao=680$~nm.
Above this region is a 75-nm-thick layer of a dielectric material with relative permittivity $\epsrel = 3.33 + 0.016i$, which acts as an anti-reflection coating. The top layer is again air-like, with relative permittivity $\epsrel = 1+10^{-9}i$.  

\begin{figure}[ht]
\begin{center}
	\begin{tabular}{cc}
		(a)\includegraphics[width=0.45\textwidth]{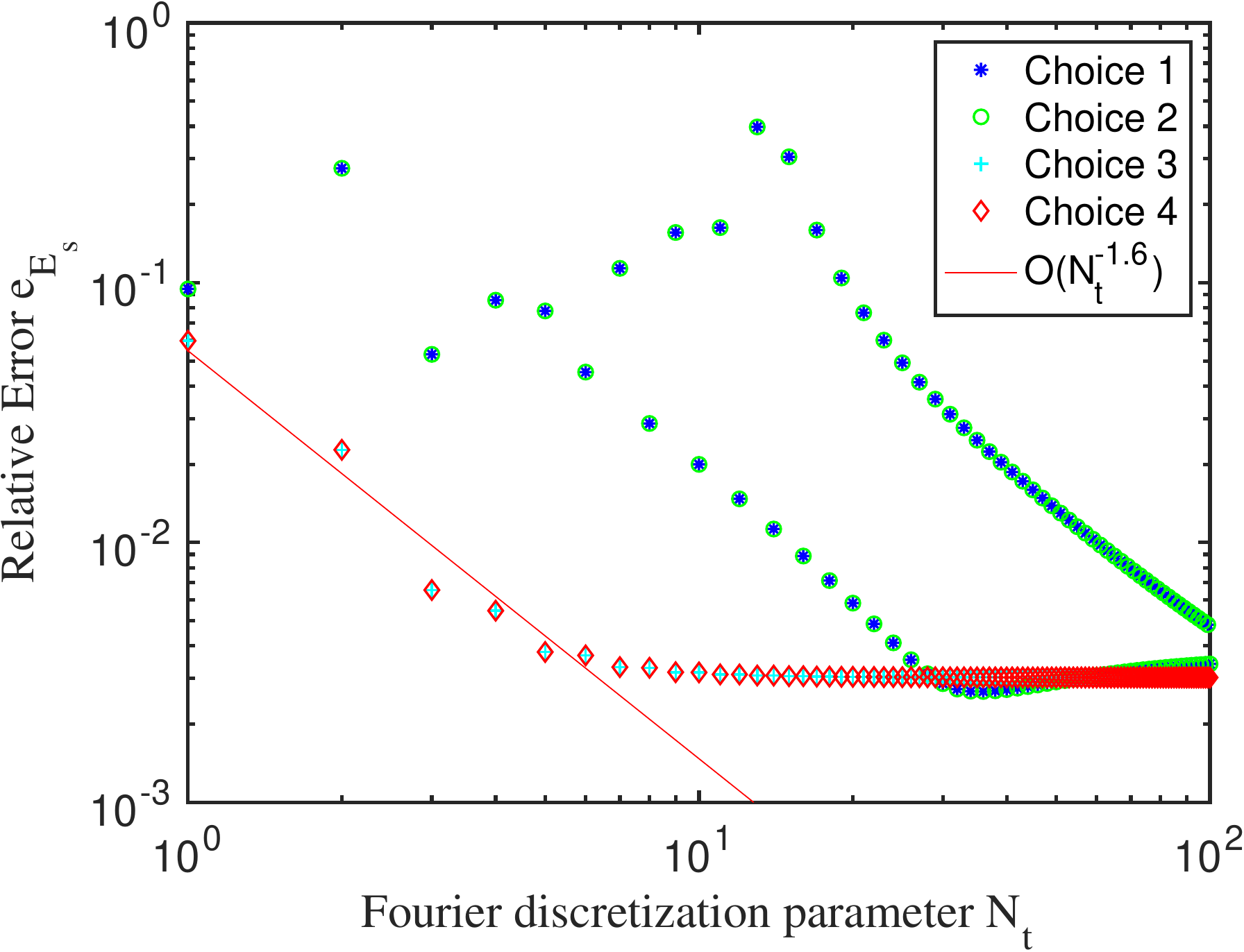}&
		(b)\includegraphics[width=0.45\textwidth]{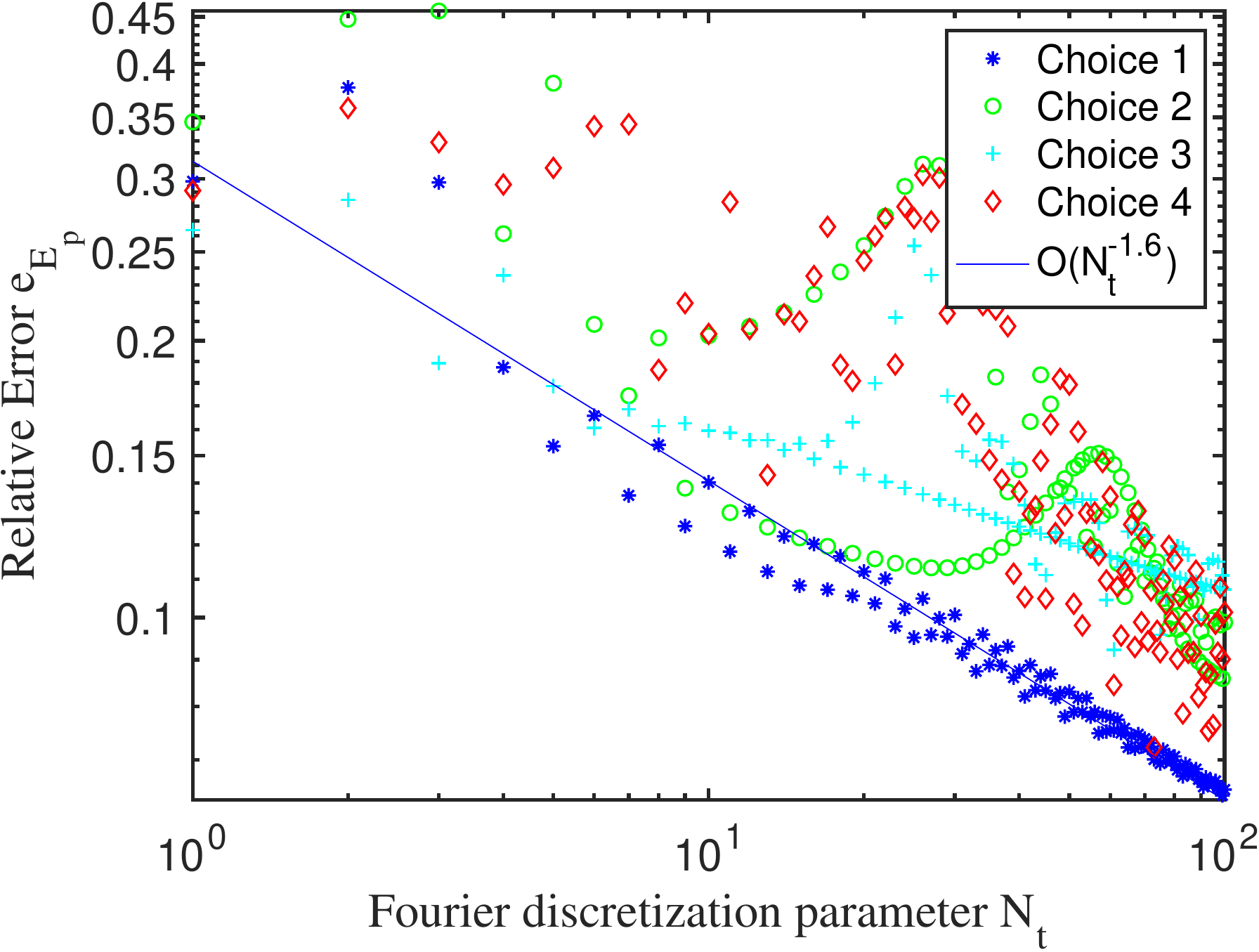}\\
	\end{tabular}
    \end{center}
 \caption{Relative errors (a) $e_{E_s}$ and (b) $e_{E_p}$ vs. $\Nt$ in the semiconductor
 region of the solar cell depicted in Fig.~\ref{Ben1}(a).
Choice 1: $\*M = \*{B}^{-1}$ and $\*N = \*{E}^{-1}$; 
Choice 2:  $\*M = \*{B}^{-1}$ and $\*N = \*{B}$; 
Choice 3: $\*M = \*{E}$ and $\*N= \*{E}^{-1}$; 
Choice 4:  $\*M = \*{E}$ and $\*N= \*{B}$.  A straight line shows the empirical order of convergence.\label{Fig:RCWAConv_p-i-n}}
\end{figure}

For this study, we used the definition
\begin{align}
\lVert f \rVert_2^2 = \int_{x=-\Lx/2}^{\Lx/2}\int_{z=0}^{\Lz}
\vert f(x,y)\vert^2\,\mathrm{d}z\,\mathrm{d}x
\label{def-fnorm-redef}
\end{align}
in Eqs.~(\ref{errEs}) and (\ref{errEp}). Four different choices of $\left\{\*M,\*N\right\}$, as described after Eq.~(\ref{eq17}),
were used in four different studies, whose results are presented
in Fig.~\ref{Fig:RCWAConv_p-i-n}.

Convergence rates for the second study are much slower than in the first study.
The best observed rates are $O(\Nt^{-1.6})$  for $s$-polarized illumination
and $O(\Nt^{-0.35})$ for $p$-polarized illumination. The slower convergence is
due to the strong singularities in the solution  near the corners of the grating. Those singularities
lead to slower convergence of the Fourier series. It should be noted that these same singularities are also
observed in the adaptive FEM solution.

Figure~\ref{Fig:RCWAConv_p-i-n}(a) contains a split in $e_{E_s}$ between successive odd and even values of
$\Nt$, when $\*M=\*E$ and $\*N=\*E^{-1}$ (Choice 3).
This split is also present, although less extreme, in an earlier work~\cite{Lalanne96}.
It is also present for $p$-polarized illumination but dramatically reduced, especially when $\*M=\*E$ and $\*N=\*B$ (Choice 1).

Lalanne \cite{Lalanne97} has already discussed in detail the complex issues behind RCWA convergence.
Our contribution is to suggest that stairstepping is not a major cause of error in RCWA, but 
it is necessary to examine convergence for different grating morphologies.  Also, RCWA is reliably convergent for $s$-polarized illumination of solar cells, and provides  slower convergence for $p$-polarized illumination.
We agree with previous studies \cite{SolanoFEM,shuba15,lokar18} that RCWA can provide the accuracy needed for photonics simulations, but our studies also
show that the convergence rates  can be quite low.   With an appropriate choice of parameters we can achieve almost 10\% error with $\Nt=10$, which became   our choice for the   optimization study reported in  Section~\ref{DEA}.

\section{Electronic Step}\label{Emodel}
Taking into account the size of the features in the solar cell, and under standard assumptions on the electron transport, such as the device being in quasi-thermal equilibrium (QTE), the transport of electrons and holes  in semiconductors can be modeled using the drift-diffusion equations~\cite{Jenny_book, Fonash_book, Brezzi2002}.
To speed up the calculations, we neglected variations with respect to $x$
and therefore used the 1D electron-hole-pair generation rate
\begin{equation}
G(z)=\frac{1}{\Lx}
\int_{-\Lx/2}^{\Lx/2}
G(x,z)\,\mathrm{d}x.\label{avgen}
\end{equation}

\subsection{Drift-Diffusion Electronic Model}\label{DD}
The DDE model comprises the following three differential equations holding
in the semiconductor region {$\Omega_{\rm el}=\left\{z\;\big\vert\;0<z<\Lz\right\}$}:
\begin{eqnarray}
&&
\frac{d}{dz}\Jn(z)= -\qe\left[G(z) - R(n,p;z)\right]\,,
\label{Eqn:dJn}
\\[5pt]
&&
\frac{d}{dz}\Jp(z)= \qe\left[G(z) - R(n,p;z)\right]\,,
\label{Eqn:dJp}
\\[5pt]
&&
 \epso\frac{d}{dz}\left[  \varepsilon_{dc} (z)\,\frac{d}{dz}\phi(z)\right]
= -\qe\left[\Nf(z) + p(z) - n(z)\right]\,.
\label{Eqn:Poisson}
\end{eqnarray}
In these equations, the gradients of the
electron current density $J_{n}(z)$, hole current density $\Jp(z)$, and electric potential $\phi(z)$
are related to the electron density $n(z)$,  hole density $p(z)$, fixed charge (charged traps and doping)  density $\Nf(z)$, electron-hole-pair generation
rate $G(z)$, and electron-hole-pair recombination
rate $R(n,p;z)$.   The zero-frequency (dc) relative permittivity is denoted by  $\varepsilon_{dc}(z)$.

The two current densities are defined as
\begin{equation}
\Jn(z) = \qe\mun(z) n(z)  \frac{d}{dz}\sfEFn(z)
\label{Eqn:Jn}
\end{equation}
and
\begin{equation}
 \Jp(z) = \qe\mup(z) p(z)  \frac{d}{dz}\sfEFp(z)\,,
 \label{Eqn:Jp} 
\end{equation}
where $\mun(z)$ and $\mup(z)$ are the electron mobility and hole mobility,
respectively; {
\begin{align}
\sfEFn(z) &= \sfEc(z) + \Vth \ln\left[\frac{n(z)}{\Nc(z)}\right]\label{Eqn:EFn}\\
\end{align}
is the electron  quasi-Fermi level; and
\begin{align}
\sfEFp(z) &= \sfEv(z) - \Vth\ln\left[\frac{p(z)}{\Nv(z)}\right]\,\label{Eqn:EFp}
\end{align}
is the   hole quasi-Fermi level. In the expressions for the quasi-Fermi levels,
$\Vth = \kB T/\qe$ is the thermal voltage of the electrons and holes with
$T$ as the temperature and  $\kB = 1.380649 \times 10^{-23}$~J~K$^{-1}$ as the Boltzmann constant,
$\Nc(z)$ is the conduction-band density of states, $\Nv(z)$
is the valence-band density of states, $\sfEc(z)$ is the conduction-band energy, and
$\sfEv(z)$ is the valence-band energy. The bandgap is defined as
\begin{equation}
\sfEg(z)=\sfEc(z)-\sfEv(z)\,
\label{Eg-def}
\end{equation}
and the intrinsic energy as}
\begin{align}
\sfEi(z) &= {\displaystyle{\frac{1}{2}\left\{\left[\sfEc(z) + \sfEv(z)\right] + \Vth
 \ln\left[\frac{\Nv(z)}{\Nc(z)}\right]\right\}}}\,.
 \label{Eqn:Ei}
\end{align}

In the Boltzmann approximation \cite{Jenny_book} the quasi-Fermi levels are
assumed to lie far from the edges of the conduction and valence bands. {Accordingly,
\begin{equation}
n(z) = \bar{n}(z) \exp\left[\frac{\sfEFn(z)- \sfEi(z)}{\Vth}\right]
\label{Eqn:ndef}
\end{equation}
and
\begin{equation}
p(z) = \bar{n}(z) \exp\left[-\frac{\sfEFp(z) - \sfEi(z)}{\Vth}\right]\,.
\label{Eqn:pdef}
\end{equation}
Here, the intrinsic carrier density is given by
\begin{equation}
\bar{n}(z) = \sqrt{\Nc(z)\Nv(z) \exp\left[\frac{\sfEg(z)}{\Vth}\right]}\,
\label{Eqn:MassAction}
\end{equation}
and the conduction-band energy is given by
\begin{align}
\sfEc(z) &= \sfE_0 - \phi(z) - \chi(z)\,,
\label{Eqn:Ec}
\end{align}
where $\chi(z)$ is the electron affinity of the semiconductor, and the (arbitrary) reference energy level $\sfE_0 $  
is often chosen as \cite{Jenny_book}
\begin{equation}
\sfE_0 = \sfEc(0) + \phi(0) + \chi(0)\,.
\end{equation}

Equations (\ref{Eqn:Jn}) and (\ref{Eqn:Jp}) may now be simplified to
\begin{equation}
\Jn(z)	= - \qe \mun(z) \left\{n(z) \frac{d}{dz}\left[\phi(z) + \phin(z)\right] - \Vth \frac{d}{dz}n(z)\right\}
\label{Eqn:Jn_Boltzmann}
\end{equation}
and
\begin{equation}
 \Jp(z)  = - \qe \mup(z) \left\{p(z) \frac{d}{dz}\left[\phi(z) + \phip(z)\right] + \Vth \frac{d}{dz}p(z)\right\}\label{Eqn:Jp_Boltzmann}\,,
 \end{equation}
 respectively.
Here, 
\begin{equation}
\label{phin-def}
\phin(z) = \chi(z) + \Vth\ln \left[\frac{\Nc(z)}{N_0}\right]
\end{equation}
and
\begin{equation}
\label{phip-def}
\phip (z)= \chi(z) + \sfEg(z) - \Vth\ln \left[\frac{\Nv(z)}{N_0}\right]
\end{equation}
 are the built-in potentials for the electrons and holes (due to variations in the material properties), respectively. Both built-in potentials as well as the electron affinity may be discontinuous with respect to position at heterojunctions in the semiconductor region.  
The baseline number density $N_0$ is arbitrary because potentials are defined up to a constant.

Equations (\ref{Eqn:dJn}), (\ref{Eqn:dJp}), (\ref{Eqn:Poisson}),
(\ref{Eqn:Jn_Boltzmann}), and (\ref{Eqn:Jp_Boltzmann}), have to be solved  
concurrently for $z\in(0,\Lz)$, in conjunction with a set of boundary conditions 
for $n(z)$, $p(z)$, and $\phi(z)$. We opted for the Dirichlet choice
\begin{align}
n(0) &= n_{0}(0),\quad &n(\Lz)&=n_{0}(\Lz),\label{Eqn:n_bc}\\
p(0) &= p_{0}(0),\quad &p(\Lz)&=p_{0}(\Lz),\label{Eqn:p_bc}\\
\phi(0) &= \phi_{0}(0),\quad &\phi(\Lz)& = \phi_{0}(\Lz) + \Vext,
\label{Eqn:Poisson_bc}
\end{align}
because it models an ideal ohmic contact \cite{Jenny_book, Fonash_book}.
 Herein, the functions $n_{0}(z)$, $p_{0}(z)$ and $\phi_{0}(z)$ are the electron density, hole density and potential at LQTE, as  discussed in Section~\ref{LQTE}, while $\Vext$ is the externally applied voltage across the terminals at the boundaries of the semiconductor region. Solution of the system of five equations enables the calculation of the current density 
\begin{equation}
\label{J-def}
J = \Jn(z) + \Jp(z)
\end{equation}
flowing \textit{uniformly} through the solar cell. 

The current density $J$   depends on the choice of $\Vext$, i.e.,
 $J\equiv J(\Vext)$. Repeating calculations for various values of $\Vext$ produces the $JV$-curve, with the maximum value of the  power density
\begin{equation}
P(\Vext) = J(\Vext)\Vext
\end{equation}
 indicating the maximum   power density 
\begin{equation}
\Pmax=\max_{\Vext} P(\Vext)
\end{equation}
obtainable from the solar cell.
This in turn gives the efficiency of the solar cell as 
\begin{align}
\eta &= \frac{\Pmax}{\Pin}\label{etadef}
\end{align}
where  $\Pin = 1000$~W~m$^{-2}$ is the incident solar power density. It is the efficiency of the solar cell that we wish to optimize.

\subsection{Local Quasi-Thermal Equilibrium}\label{LQTE}
If a region containing a homogeneous semiconductor is charge-free and   isolated from any external influences (e.g., $G \equiv 0$ and $\Vext = 0$), the distributions of electrons and holes with respect to energy tend to the Fermi distribution, with the Fermi-level lying at the intrinsic energy level, i.e., $\sfEFn = \sfEFp = \sfE_F$. Then, $\Jn=\Jp$ are identically zero and the semiconductor is said to be in LQTE. This  concept is used to model the ideal ohmic boundary conditions mentioned in Section~\ref{DD}.
	
With $n_0(z)$ and $p_0(z)$ denoting the  electron and hole densities, respectively, at LQTE, Eqs.~(\ref{Eqn:ndef})--(\ref{Eqn:MassAction}) can be 
combined to give the \textit{mass-action} equation
\begin{align}
\bar{n}^2(z) = n_0(z) p_0(z)\,;\label{Eqn:massaction_n0p0}
\end{align}
furthermore,
\begin{align}
\Nf(z) + p_0(z)-n_0(z) &= 0\,,
\label{Eqn:rho0}
\end{align}
as the semiconductor is charge-free.
Combining Eqs.~(\ref{Eqn:massaction_n0p0}) and (\ref{Eqn:rho0}), we get
\begin{align}
n_0(z) &= \frac{\Nf(z) + \sqrt{{\Nf(z)^2}+4\bar{n}(z)^2}}{2}\label{Eqn:n0}\\
\intertext{and}
p_0(z) &= \frac{-\Nf(z) + \sqrt{{\Nf(z)^2}+4\bar{n}(z)^2}}{2}\,\label{Eqn:p0}.
\end{align}

At {LQTE}, the electron and hole quasi-Fermi levels coincide and are {uniform}. Thus, {
\begin{align}
\sfEFn(\upsilon) &= \sfEc(\upsilon) + \Vth\ln\left[\frac{n_0(\upsilon)}{\Nc(\upsilon)}\right] \\
\intertext{and}
\sfEFp(\zeta) &= \sfEv(\zeta) - \Vth\ln\left[\frac{p_0(\zeta)}{\Nv(\zeta)}\right],
\end{align}
must be equal to each other}
for all $\upsilon  \in (0,\Lz)$ and $\zeta \in (0,\Lz)$. 
Substitution of Eqs.~(\ref{Eg-def}) and (\ref{Eqn:Ec}) into $\sfEFn(\upsilon)=\sfEFp(\zeta)$ yields {
\begin{align}
-\phi_0(\upsilon) - \chi(\upsilon) + &\Vth\ln\left[\frac{n_0(\upsilon)}{\Nc(\upsilon)}\right] =\nonumber\\
 & -\phi_0(\zeta) - \chi(\zeta) -\sfEg(\zeta) + \Vth\ln\left[\frac{p_0(\zeta)}{\Nv(\zeta)}\right].
 \label{eq58-al}
\end{align}
Equation~(\ref{eq58-al}) simplifies to
\begin{align}
&- \chi(0) + \Vth\ln\left[\frac{n_0(0)}{\Nc(0)}\right] =
  -\phi_0(z) - \chi(z) -\sfEg(z) + \Vth\ln\left[\frac{p_0(z)}{\Nv(z)}\right]
 \label{eq59-al}
\end{align}
with the choice $\upsilon= 0$ and $\zeta = z$, and to
\begin{align}
-\phi_0(z) - \chi(z) + &\Vth\ln\left[\frac{n_0(z)}{\Nc(z)}\right] =
  - \chi(0) -\sfEg(0) + \Vth\ln\left[\frac{p_0(0)}{\Nv(0)}\right]
 \label{eq60-al}
\end{align}
with the choice $\upsilon = z$ and $\zeta = 0$,   where we have chosen   $\phi_0(0)=0$ without loss of generality.
Subtracting Eq.~(\ref{eq59-al}) from Eq.~(\ref{eq60-al}), we get}
\begin{align}
 \phi_0(z)  &=  - [\chi(z) - \chi(0)] - \frac{1}{2}[\sfEg(z) - \sfEg(0)] \nonumber\\
 & +\frac{1}{2}\Vth\ln\left[\frac{n_0(z)}{n_0(0)}\right] + \frac{1}{2}\Vth\ln\left[\frac{p_0(z)}{p_0(0)}\right] \nonumber\\
 & -\frac{1}{2}\Vth\ln\left[\frac{\Nc(z)}{\Nc(0)}\right] - \frac{1}{2}\Vth\ln\left[\frac{\Nv(z)}{\Nv(0)}\right]\,.\label{Eqn:phi0}
\end{align}
 
 Equations (\ref{Eqn:n0}), (\ref{Eqn:p0}), and (\ref{Eqn:phi0}) give the initial functions
needed to specify the boundary data  identified in Eqs.~(\ref{Eqn:n_bc})--(\ref{Eqn:Poisson_bc}).   However, in view of the LQTE assumption (\ref{Eqn:rho0}),   these initial functions do not satisfy the system of ODEs (\ref{Eqn:dJn})--(\ref{Eqn:Poisson}). 

\subsection{Recombination}\label{Recom}
Recombination of an electron and a hole can take place through several different
mechanisms \cite{Jenny_book,Fonash_book}. The electronic step in our program
accommodates three different contributions to $R(n,p;z)$.
 
\subsubsection{Radiative Recombination}
Radiative recombination occurs when an electrons and a hole recombine across the full bandgap, releasing the energy as a photon with energy equal to the bandgap. At QTE, radiative recombination is identically balanced by electrons being thermally excited across the bandgap. The (net) radiative recombination rate is given by
\begin{align}
R_{rad}(n,p;z) = \frac{\alpha_{rad}(z)}{\bar{n}^2(z)} \left[n(z)p(z) - \bar{n}^2(z)\right]\,,
\label{rad}
\end{align}
where $\alpha_{rad}(z)$   depends on the semiconductor material. It should be noted that {$R_{rad}(z)\equiv 0\, \forall z\in\left[0,\Lz\right]$} at LQTE.

\subsubsection{Shockley--Read--Hall (SRH) Recombination}
The SRH recombination contribution is due  to   electrons and holes recombining via an intermediate gap state. It is modeled by
 \begin{align}
R_{SRH}(n,p;z) = \frac{n(z)p(z) - \bar{n}^2(z)}{\taup(z)\left[n(z) + n_1(z)\right]+ \taun(z)\left[p(z)+p_1(z)\right]}\,,\label{SRH}
\end{align}
where $n_1(z)$ and $p_1(z)$ are, respectively, the electron and hole densities  at the trap energy level. If this level is the intrinsic energy level $\sfEi(z)$, then $n_1(z) = p_1(z) = \bar{n}(z)$
from Eqs.~(\ref{Eqn:ndef}) and (\ref{Eqn:pdef}).  The functions
$\taun(z)$ and $\taup(z)$ are material dependent.

\subsubsection{Auger Recombination}
The Auger recombination contribution arises from a three-particle recombination pathway, occurring when an electron and
a hole recombine across the bandgap, with the released energy transferred to a third charge carrier which is excited away from both band edges. This recombination rate is given by 
 \begin{align}
R_{Aug}(n,p;z) &= C_{\rm n}(z)n(z)\left[n(z)p(z) - \bar{n}^2(z)\right] 
\nonumber
\\
&+ C_{\rm p}(z)p(z)\left[n(z)p(z) - \bar{n}^2(z)\right]\,,\label{Aug}
\end{align}
where the functions $C_{\rm n}(z)$ and $C_{\rm p}(z)$ are material dependent.

All three contributions add, so that
\begin{equation}
R (n,p;z) = R_{rad}(n,p;z)+ R_{SRH}(n,p;z)+ R_{Aug}(n,p;z)\,.
\end{equation}
It can be seen that the total recombination $R = 0$ at LQTE irrespective of the choice of the material parameters. For all results reported in this paper, we set {$\alpha_{rad}(z)\equiv0$ and $C_{\rm n}(z)=C_{\rm p}(z)\equiv0\, \forall z \in\Omega_{\rm el}$} and only the SRH recombination was activated,
in order to present illustrative results without the expenditure of significant computational time.
 
\subsection{Heterojunctions: Continuous Quasi-Fermi Levels}\label{Sec:Hetero}
When there is a jump in either  $\chi(z)$ or $\sfEg(z)$ or both,
and thus in $\phin(z)$ or $\phip(z)$, then $n(z)$ and  $p(z)$ may also
have   jumps. A heterojunction is formed at the discontinuity \cite{Hetero}. A commonly used way to model this discontinuity requires the assumption of continuous quasi-Fermi levels,
whereas another commonly used way requires consideration of thermionic emission at the discontinuity \cite{Therm}. We chose to implement the continuous quasi-Fermi level (CQFL) model, although the thermionic-emission model
is also compatible with the HDG method.

The CQFL model uses the limit of a continuum model to quantify the jump at a heterojunction. To understand the resulting jump conditions, let us examine a small region $0< z \leq \delta$ in which the electron affinity $\chi(z)$ changes linearly by $\Delta \chi$, while $\Nc(z)$  is uniform. Then, differentiating Eq.~(\ref{phin-def}) with respect to $z$ yields
\begin{align}
\label{eq58-AL}
\frac{d}{dz}\phin(z)= \frac{\Delta \chi}{\delta}\,
\end{align}
in this region.
Furthermore, we assume that $G \delta \approx R \delta \approx 0$ in this region so that
Eq.~(\ref{Eqn:dJn}) simplifies to
\begin{align}
\frac{d}{dz}\Jn(z)&= 0\,;
\end{align}
hence, {$\Jn(z)=\Jn\zero$.
Finally, we also assume 
that both $\mun(z) \approx \mun\zero$ and $\phi(z) \approx \phi\zero$} are uniform in this region,
so that Eq.~(\ref{Eqn:Jn_Boltzmann}) simplifies to
\begin{equation}
\Jn\zero	= - \qe \mun\zero   n(z) \frac{\Delta \chi}{\delta} + \qe \Vth \mun\zero   \frac{d}{dz}n(z) \,,
\label{eq60-AL}
\end{equation}
after using Eq.~(\ref{eq58-AL}).

The solution of Eq.~(\ref{eq60-AL}) is
\begin{align}
\label{eq69-al}
n(z) 	&= C \exp\left(\frac{\Delta \chi}{\Vth}\frac{z}{\delta}\right) - \frac{\Jn\zero \delta}{\qe\mun\zero \Delta \chi},
\end{align}
where $C$ is the constant of integration. As $n(z) = \nL$ at  $z=0$, we get
\begin{align}
\label{eq70-al}
	n (z)	&= \left(\nL + \frac{\Jn\zero \delta}{\qe \mun\zero \Delta \chi}\right) \exp\left(\frac{\Delta \chi}{\Vth} \frac{z}{\delta}\right) - \frac{\Jn\zero \delta}{\qe\mun\zero \Delta \chi}
\end{align}
from Eq.~(\ref{eq69-al}).
As $n(z) = \nR$ at  $z=\delta$, Eq.~(\ref{eq70-al}) yields
\begin{align}
\nR 	
	& = \left(\nL + \frac{\Jn\zero\delta}{\qe \mun\zero \Delta \chi}\right) \exp\left[\frac{\Delta \chi}{\Vth}\right] - \frac{\Jn\zero \delta}{\qe \mun\zero \Delta \chi}.
	\label{eq71-al}
\end{align}

In the limit of $\delta \to 0$, we   get the jump condition
\begin{align}
	\nR = \exp\left[\frac{\Delta \chi}{\Vth}\right] \nL
	\label{hn}
\end{align}
from Eq.~(\ref{eq71-al}). The same analysis can be performed for the hole density $p(z)$, which results in the jump condition
\begin{align}
\pR = \exp\left[-\frac{\Delta \chi + \Delta \sfEg}{\Vth}\right] \pL,\label{hp}
\end{align}
where {$\Delta \sfEg$ is the change in the bandgap over the region with thickness $\delta$.
Equations (\ref{hn}) and (\ref{hp}) need to be enforced across any junction {over} which $\chi$ or $\sfEg$ jump. This is accomplished using the HDG method described in Section~\ref{HDG}.

\subsection{Dimensionless Formulation}\label{dimless}
{The equations of the DDE model} are scaled in order to ensure that the internal variables are dimensionless \cite{Brezzi2002}, the scaling parameters being provided in Table~\ref{Tab:Data}. Equations (\ref{Eqn:dJn}) and (\ref{Eqn:dJp}) are thus non-dimensionalized to\footnote{At the risk of confusion, we have not used new notation for the dimensionless variables in order to keep the text clean.}
\begin{align}
\frac{d}{dz}\Jn(z)&= -G(z) + R(n,p;z) \label{Jnp}\\
\intertext{and}
\frac{d}{dz}\Jp(z)&= G(z) - R(n,p;z), \label{Jpp}
\end{align}
respectively, and Eqs. (\ref{Eqn:Jn_Boltzmann}) and (\ref{Eqn:Jp_Boltzmann})
to
\begin{align}
\Jn(z) & = -\mun(z) \left\{n(z) \frac{d}{dz}[\phi(z) + \phin(z)] 
-\frac{d}{dz}n(z)\right\}\label{Jn}\\
\intertext{and}
\Jp(z) & = -\mup(z) \left\{p(z) \frac{d}{dz}[\phi(z)+  \phip(z)]
+ \frac{d}{dz}p(z)\right\},\label{Jp}
\end{align}
respectively. Electing to keep $\varepsilon_{dc}(z)\equiv \varepsilon_{dc}^0\,\forall z \in\left[0,\Lz\right]$ uniform in the semiconductor region, 
we obtain the non-dimensionalized form of Eq.~(\ref{Eqn:Poisson}) as
\begin{equation}
-\lambda^2\frac{d^2}{dz^2}\phi(z)= \Nf(z)+p(z)-n(z)\,,
\label{eq71-AL}
\end{equation}
where the Poisson constant
\begin{equation}
\lambda^2  =  \frac{\epso\varepsilon_{dc}^0 \phis}{\qe \Ns \Ls^2}\,
\end{equation}
{is dimensionless.}
We define {the dimensionless d.c. electric field}
\begin{equation}
\Edc (z) = -\lambda^2\frac{d}{dz}\phi(z)
\label{E}
\end{equation}
and transform Eq.~(\ref{eq71-AL}) to
\begin{equation}
\frac{d}{dz}\Edc(z)= \Nf(z)+p(z)-n(z)\,,
\label{Ep}
\end{equation}
which has the same form as Eqs.~(\ref{Jnp}) and (\ref{Jpp}). These three ODEs
are  discretized using the HDG method, as described in Section~\ref{HDG}.

\begin{table}[ht]
	\begin{center}	
		\begin{tabular}{|l|c|c|c|}\hline
			Scaling 	& Symbol & Representative Value & Variable(s)
			\\\hline\hline
			Length & $\Ls$ & $\mathcal{O}(\Lz) \approx 2\times10^{-4}$~cm
			& $z$\\
\hline
			Mobility & $\mus$ & $1.5$~cm$^2$~V$^{-1}$~s$^{-1}$ 
			& $\mu_{\rm n,p}$\\
\hline
			Voltage & {$\phis=\Vth$} & $0.0259$~V
			& $\phi$, $\phi_{\rm n,p}$, $\chi$, $\Vext$, \\
						  & &  & $\sfE_{\rm i,c,F_{\rm n},F_{\rm p},g,v}$\\
\hline
			Density & $\Ns$& $\max(\vert\Nf\vert) \approx 10^{16}$~cm$^{-3}$ 
			& $n$, $p$, $\bar{n}$, \\
			  & &  & $n_{0,1}$, $p_{0,1}$,{$N_{\rm c,f,v}$}\\
			\hline
			Poisson const. & $\lambda^2$ & $4.1717\times10^{-4}$& --\\
			\hline
			Density rate & {$\Gs=\dfrac{\mus\Ns\phis}{\Ls^2}$} & $9.6945\times 10^{21}$ cm$^{-3}$~s$^{-1}$ & $G$, $R$\\
			\hline
			Current & {$\Js=\dfrac{\qe\mus\Ns\phis}{\Ls}$} & $0.3106$~mA~cm$^{-2}$
			& $J$, $J_{\rm n,p}$\\ 
			\hline
		\end{tabular}
		
		\caption{Table of  scaling parameters and applicable variables.} \label{Tab:Data}
	\end{center}
\end{table}

\section{Hybridizable Discontinuous Galerkin Formulation}\label{HDG}
{The hybridizable discontinuous Galerkin (HDG) method~\cite{FuQiuHDG,chen03}  possesses several features which are advantageous for the implementation of the DDE model for solar cells~\cite{Brinkman}.} The relaxation of the requirement that the solution is continuous at the boundaries of elements permits solutions where the variables have strong gradients and second derivatives. It also naturally allows discontinuities in the solution due to discontinuous material parameters, such as those which occur at heterojunctions discussed in Section~\ref{Sec:Hetero}.

\subsection{Generalized Transport System}\label{sec:gts}
To simplify the description of the numerical scheme used to solve  Eqs.~(\ref{Jnp})--(\ref{Jp}), (\ref{E}), and (\ref{Ep}), it is useful to partition them into three sets, each with the form
\begin{align}
\sfc_1 \sfJ & = \sfc_2 \sfm \frac{d}{dz}(\sfZ+ \sfvm) +\frac{d}{dz} \sfm, &\label{gt1}\\
\frac{d}{dz}\sfJ &= \sff(n,p)  + \sfc_3(p-n),&\label{gt2}
\end{align}
along the associated Dirichlet boundary conditions
\begin{equation}
\left.
\begin{array}{l}
\sfm(0) = \sfm_0(0)\\[5pt]
\sfm(\Lz)=\sfm_0(\Lz)
\end{array}\right\}\,.\label{gt3}
\end{equation}
Here, $\sfZ$ is either given by the {solution of another set   or vanishes,} $\sfvm$ is a given function of position, and $\sff(n,p)$ is also given.  To get each of the three sets of equations in the dimensionless DDE model provided in Section~\ref{dimless},
we choose the parameters and functions as in Table~\ref{tab_equiv}.

\begin{table}
\begin{tabular}{||c|ccc||}\hline\hline
&Eqs. (\ref{Jnp}) and (\ref{Jn}) &Eqs. (\ref{Jpp}) and (\ref{Jp})
&Eqs. (\ref{E}) and (\ref{Ep})\\\hline
$\sfJ$ & $\Jn$&$\Jp$& $E$\\
$\sfm$ & $n$&$p$&$\phi$\\
$\sfZ$ & $\phi$ &  $\phi$ & 0\\
$\sfvm$&$\phin$& $\phip$&0\\
$\sfc_1$ & $1/\mun$ & $-1/\mup$ & $-1/\lambda^2$\\
$\sfc_2$ & $-1$&$1$&$0$\\
$\sfc_3$&0&0&1\\
\hline\hline
\end{tabular}
\caption{Correspondence between the generalized transport system of Eqs.~(\ref{gt1})--(\ref{gt3})
and Eqs.~(\ref{Jnp})--(\ref{Ep}) of the dimensionless DDE model.}
\label{tab_equiv}
\end{table}

Henceforth, we refer to {Eqs.~(\ref{gt1})--(\ref{gt3})} as the \textit{generalized transport system}. 
We now proceed to design a numerical scheme to discretize this system of equations. 
 
\subsection{Discretization of Electronic Domain}\label{HDGD}

The electronic domain {$\Omega_{\rm el}$} is covered by a mesh $\mathcal{S} = \{\mathcal{N}, \mathcal{T}\}$ of $\Nz+1$ nodes $\mathcal{N}$ and $\Nz$ elements $\mathcal{T}$.
The element $s_\gamma=({z}_{\gamma-1} ,{z}_\gamma)$ lies between the nodes ${z}_{\gamma-1}$
and ${z}_{\gamma}$, 
\begin{align}
\mathcal{N} = \left\{ {z}_{\gamma} |{z}_{\gamma}\in\Omega_{\rm el}, \gamma \in\{0,..,\Nz\}\right\}.
\end{align}
is the set of all nodes,
and
\begin{align}
\mathcal{T} = \left\{ s_\gamma | s_\gamma=({z}_{\gamma-1} ,{z}_\gamma) \subset \Omega_{\rm el},
\gamma\in\{1,...,\Nz\}\right\}
\end{align}
is the set of all elements.
The generalized transport system is discretized using the space
\begin{equation}
\mathbb{V}_h=\left\{ v_h\in L^2({\Omega_{\rm el}})\;\big\vert\; v_h\vert_{s_\gamma}\in \mathbb{P}_{\Pdeg}
~\forall s_\gamma \in \mathcal{T}\right\}\,,
\end{equation}
where $\mathbb{P}_{\Pdeg}$ is the set of polynomials  of degree
$\Pdeg$ in $z$ {(and so that $\mathbb{V}_h$ is a space of discontinuous piecewise polynomials)}, 
\begin{equation}
 \mathbb{W}_h = \{w_h \in w_h|_{z_\gamma} \in \mathbb{R}^+, ~\forall {z}_{\gamma} \in \mathcal{N} \},
\end{equation}
and $\mathbb{R}^+$ is the set of positive real numbers.
Although a different  value of $\Pdeg$ may be chosen for each  element $s_{\gamma}$, for simplicity of notation we choose the same value of $\Pdeg$ for all elements.

We  seek the numerical solutions $(\sfJ_h, \hat{\sfJ}_h, \sfm_h, \hat{\sfm}_h)\in \mathbb{V}_h\times \mathbb{W}_h \times \mathbb{V}_h \times \mathbb{W}_h$ {where $\hat{\sfJ}_h$ and $\hat{\sfm}_h$ are the vectors of numerical approximations to $\sfJ$ and $m$ at the nodes in $\cal N$.}
Multiplying Eq.~(\ref{gt1}) by  a test function $\psi \in\mathbb{V}_h$ and Eq.~(\ref{gt2}) by  
a test function $\eta \in\mathbb{V}_h$, and then integrating over ${\Omega_{\rm el}}$, we obtain
\begin{align}
\sum_{\gamma=1}^{\Nz}\pt{\sfc_1 \sfJ,\psi} & = \sfc_2\sum_{\gamma=1}^{\Nz} \pt{\sfm\frac{d}{dz}\sfZ,\psi} + \sfc_2 \sum_{\gamma=1}^{\Nz} \pt{\sfm\frac{d}{dz}\sfvm,\psi}\nonumber\\&\qquad + \sum_{\gamma=1}^{\Nz}\pt{\frac{d}{dz}\sfm,\psi} &\label{Eqn:FirstTestJ}\\
\intertext{and}
\sum_{\gamma=1}^{\Nz}\pt{\frac{d}{dz}\sfJ,\eta}&= \sum_{\gamma=1}^{\Nz}\pt{\sff,\eta}  + \sfc_3\sum_{\gamma=1}^{\Nz}\pt{(p-n),\eta},\label{Eqn:FirstTestn}
\end{align}
where the notation
\begin{equation}
\pt{\eta,\pi} = \int_{s_\gamma} \eta(z) \psi(z)\,{\rm{}d}z.
\end{equation}

For functions $\eta \in \mathbb{V}_h$, $\psi \in \mathbb{V}_h$, and $\hat{\eta} \in \mathbb{W}_h$,  we define 
\begin{equation}
\pt{\frac{d}{dz} \eta, \psi} = \ipt{\hat{\eta},\psi}  - \pt{\eta, \frac{d}{dz}\psi}
\end{equation}
via integration by parts,
with
\begin{equation}
\ipt{\hat{\eta},\psi}  = \hat{\eta}_{\gamma} \psi|_{\gamma,R} - \hat{\eta}_{\gamma-1} \psi|_{\gamma,L}\,
\end{equation}
and
 \begin{equation}
 \left.\begin{array}{l}
 \psi\vert_{\gamma,L} = \lim_{z{\to} {z}_{\gamma-1}^+}\psi(z)
 \\[5pt]
 \psi\vert_{\gamma,R} = \lim_{z{\to} {z}_\gamma^-}\psi(z)
 \end{array}\right\}\,.
 \end{equation}
Then, Eq.~(\ref{Eqn:FirstTestJ}) can be written as
\begin{align}
{\sum_{\gamma=1}^{\Nz}\Bigg[}
\pt{\sfc_1 \sfJ,\psi} &- \sfc_2 \pt{\sfm\frac{d}{dz}\sfZ,\psi} - 
\sfc_2 \pt{\sfm\frac{d}{dz}\sfvm,\psi} + \pt{\sfm,\frac{d}{dz}\psi} \nonumber\\&- \ipt{\hat{\sfm},\psi}
{\Bigg]} =0 \\
\intertext{and Eq.~(\ref{Eqn:FirstTestn}) as}
{\sum_{\gamma=1}^{\Nz}\Bigg[}
 \pt{\sfJ,\frac{d}{dz}\eta}& + \sfc_3\pt{p,\eta} - \sfc_3\pt{n,\eta} + \pt{\sff,\eta}  - \ipt{\hat{\sfJ},\eta}
{\Bigg]}=0\,.
\end{align}

To complete the discretization, we need to enforce approximate continuity of the discrete flux.  Therefore, at each node we define the discrete flux as
\begin{align}
\hat{\sfJ}_\gamma = \sfJ({z}_\gamma) + \tau_m({z}_\gamma)[\sfm({z}_\gamma)-\hat{\sfm}_\gamma]\,,
\end{align}
where $\zeta({z}_{\gamma}) = \lim\limits_{z\to {z}^{\pm}_\gamma} \zeta(z)$ and $\tau_m({z}_\gamma)\geq 0$ is a hybridization (or penalty) function to be chosen.
Continuity of the flux requires
\begin{align}
\sfJ\vert_{\gamma+1,L} + \tau_m(\sfm-\hat{\sfm}_\gamma)\vert_{\gamma+1,L} = \sfJ\vert_{\gamma,R} + \tau_m(\sfm-\hat{\sfm})\vert_{\gamma,R}
\end{align}
at each interior node.
Finally, the boundary conditions are enforced by requiring
\begin{equation}
\hat{\sfm}=\sfm_0 \label{Eqn:BChat}
\end{equation}
at the ends of $\Omega$.

By defining the jump operator as
\begin{align}
\dsbg{\xi} = \xi\vert_{\gamma,R} - \xi\vert_{\gamma,L}\,,
\end{align}
the problem becomes that of finding a discrete solution {
${(\sfJ_h,\sfm_h,\hat{\sfm})} \in \mathbb{V}_h\times \mathbb{V}_h\times \mathbb{W}_h$ such that
\begin{align}
{\sum_{\gamma=1}^{\Nz}\Bigg[}
 \pt{\sfc_1 \sfJ_h,\psi} &
- \sfc_2\pt{\sfm_h\,\frac{d}{dz}\sfZ,\psi} - \sfc_2\pt{\sfm_h\frac{d}{dz}\sfvm,\psi} + 
\pt{\sfm_h,\frac{d}{dz}\psi} \nonumber\\
	& - \ipt{\hat{\sfm},\psi}{\Bigg]}=0,
	\label{Eqn:system2start}\\
	{\sum_{\gamma=1}^{\Nz}\Bigg[}
 \pt{\sfJ_h,\frac{d}{dz}\eta} &
	 + \sfc_3\pt{p_h,\eta} - \sfc_3\pt{n_h,\eta} + \pt{\sff,\eta} &\nonumber\\
	& -\ipt{\sfJ_h,\eta} - \ipt{\tau_m \sfm_h,\eta} + \ipt{\tau_m \hat{\sfm},\eta}{\Bigg]}=0,\label{Eqn:system2next}\\
	 \dsbg{\sfJ_h} &+ \dsbg{\tau_m \sfm_h} - \dsbg{\tau_m}\hat{\sfm}
	 =0
	 \;\mbox{ at all interior nodes}, \label{Eqn:system2nextnext}
\end{align}
for} all $(\psi,\eta)\in \mathbb{V}_h\times \mathbb{V}_h$, and with $\hat{\sfm}=\sfm_0\mbox { for }
{z\in\left\{0,\Lz\right\}}$.  Obviously, there are
$2\Nz(\Pdeg+1)$ degrees of freedom for the functions $\sfJ$ and $\sfm$, and $\Nz-1$
degrees for $\hat{\sfm}$.  In turn, the test functions {$\psi$ and $\eta$} give 
{$2\Nz(\Pdeg+1)$} equations from Eqs.~(\ref{Eqn:system2start}) and
(\ref{Eqn:system2next}, while Eqs.~(\ref{Eqn:system2nextnext}) give a further $\Nz-1$ equations. 
For a linear problem of this type, existence and uniqueness of the solution have been proven and the error estimates have been derived, for example, by Fu~\textit{et al.} \cite{FuQiuHDG}; see also 
Cockburn~\textit{et al.} \cite{CockburnHDG}. The nonlinear problem 
defined here remains to be analyzed.

To help with conditioning and also to provide a useful choice of discretization points for the nonlinear problem, we choose
the interpolation points for $\mathbb{V}_h$ on each 
element to be the $\Pdeg+1$ Gauss--Lobatto points~\cite{str73}. We now replace functions $\sfJ_h$, $\sfm_h$, $\sfZ_h$, $n_h$, and $p_h$ by their finite-element approximations with the form
\begin{align}
\label{eq100-al}
&\zeta_h{(z)} = \sum_{\ell=1}^{\Nz(\Pdeg+1)} \zeta_\ell\,\psi_\ell{(z)},
\end{align}
where $\psi_\ell \in \mathbb{V}_h$ are piecewise polynomial finite-element basis functions that are each compactly supported on one of the elements $s_{\gamma}$, and $\zeta_{\ell} \in \mathbb{R}$ are constants to be determined.
Upon substitution of these into Eqs.~(\ref{Eqn:system2start})--(\ref{Eqn:system2nextnext}), and choosing the test functions $\eta$ and $\psi$ to also be the finite-element basis functions $\psi_j \in \mathbb{V}_h$ to be piecewise polynomials that are each compactly supported on element $s_{\gamma}$, we obtain the HDG equations
\begin{align}
{\sum_{\gamma=1}^{\Nz}\Bigg[}
\sfJ_\ell\pt{\sfc_1 \psi_\ell,\psi_j} &- \sfc_2 \sfm_\ell \sfZ_k \pt{\psi_\ell \frac{d}{dz}\psi_k,\psi_j} - 
\sfc_2 \sfm_\ell\pt{\frac{d}{dz}\sfvm \psi_\ell,\psi_j} \nonumber\\&
+  \sfm_\ell\pt{\psi_\ell,\frac{d}{dz}\psi_j} - \ipt{\hat{\sfm},\psi_j}
{\Bigg]}=0, \label{Eqn:HDG_start}\\
{\sum_{\gamma=1}^{\Nz}\Bigg[}
\sfJ_\ell\pt{\psi_\ell,\frac{d}{dz}\psi_j} &
  - \sfc_3 \left[n_\ell\pt{\psi_\ell,\psi_j} + p_\ell\pt{\psi_\ell,\psi_j}\right] \nonumber
\\&  - \sfm_\ell\ipt{\tau_m \psi_\ell,\psi_j} + \pt{\sff,\psi_j} + \ipt{ \hat{\sfm},\psi_j}
{\Bigg]}=0,\\
 \sfJ_\ell \dsbg{\psi_\ell} &+ \sfm_\ell\dsbg{\tau_m \psi_\ell} - \dsbg{\tau_m}\hat{\sfm_\ell} {=0}
  \mbox{ at interior nodes},\\
\hat{\sfm}_{\gamma}&=\sfm_0\mbox{ at the end points of }\Omega_{\rm el}.
\label{Eqn:HDG_end},
\end{align}
for $1\leq j \leq \Nz(\Pdeg+1)$, {$1\leq k \leq \Nz(\Pdeg+1)$,} and $1\leq  \ell\leq \Nz(\Pdeg+1)$.
In Eqs.~(\ref{Eqn:HDG_start})--(\ref{Eqn:HDG_end})
and henceforth, we use the Einstein summation convention  (i.e., summation is implied over any duplicated index),  column vectors  denoted as {$\mfv$} comprise elements $v_k$, and matrices denoted as {$\mfM$} comprise elements $M_{k,\ell}$.

The full set of drift-diffusion equations (\ref{Jnp})--(\ref{Ep})
along with the boundary conditions (\ref{Eqn:n_bc})--(\ref{Eqn:Poisson_bc})
 is now discretized
by three copies of the above equations with function
and parameter choices from Table~\ref{tab_equiv}.
There are two usual ways to solve the resulting  nonlinear system: Gummel iteration  \cite{Brezzi2002} and Newton's method \cite{isa66}. We found that straightforward Gummel iteration did not converge reliably for the parameter values in our simulations. So we elected to use Newton's method with a scheme of the homotopy type  \cite{nocedal}  to help
provide a good initial guess.

\subsection{Nonlinear Convection \label{Sec:Nonlinear}}
{A practical difficulty with implementing the scheme arises when} dealing with the nonlinear term 
\begin{equation}
{\underline \psi_{\sfm}=\sum_{\gamma=1}^{\Nz}\Bigg[}\sfm_\ell \sfZ_k \pt{\psi_\ell \frac{d}{dz}\psi_k,\psi_j}{\Bigg]}. 
\label{108-AL}
\end{equation}
The integral  
\begin{equation}
{\Psi^{0,0,1}_{\ell,j,k}=\sum_{\gamma=1}^{\Nz}\Bigg[}\pt{\psi_\ell \frac{d}{dz}\psi_k,\psi_j}{\Bigg]}
\end{equation}
on the right side of Eq.~(\ref{108-AL})
 {is an element of
 a $[\Nz(\Pdeg+1)]^3$ tensor denoted by $\~\Psi^{0,0,1}$,} where the $0,0,1$ superscript specifies that {$\psi_k$ is differentiated once}.
 The non-linear term can then be seen as the product
\begin{align}
\underline \psi_{\sfm} \equiv \Psi^{0,0,1}_{\ell,j,k} \sfm_\ell \sfZ_k\,,
\end{align}
where the indices $\ell$, $j$,  and $k$ select the elements of the tensor.
Although manipulation of tensors with rank greater than $2$  is not natively supported in Matlab,  the tensor 
{$\~\Psi^{0,0,1}$}
is very sparse as it is block diagonal; hence, the calculation speed can be increased by storing it as a non-square matrix. 
	
As {each of the finite-element basis functions $\psi_j$ is} compactly supported on just one interval, a \textit{local} tensor 
\begin{equation}
\~\Phi^{0,0,1}_\gamma \equiv \pt{\psi_\ell \frac{d}{dz}\psi_k,\psi_j}
\end{equation} 
with dimensions $[(\Pdeg+1)]^3$ can be formed. This local tensor, which is not necessarily sparse, can be rewritten as a $[\Pdeg+1]\times [(\Pdeg+1)^2]$ rectangular matrix {$\~\Psi^{0,0,1}_\gamma$.} Finally, a global matrix {$\~\Psi^{0,0,1}$} with dimension $\Nz(\Pdeg+1)\times (\Nz(\Pdeg+1))^2$ to represent {$\~\Psi^{0,0,1}$}
 is formed by repeating the local matrix as a block diagonal.

 In order to calculate the term {$\underline \psi_\sfm$}, the vector of coefficients {$\underline\sfm_\ell$} is reshaped into a $\left[\Nz(\Pdeg+1)\right]^2 \times \Nz(\Pdeg+1)$ matrix {${\=M}$:} the block diagonal is formed from $\Pdeg+1$ blocks of {$\sfm_\ell$}, with each block repeated $\Nz$ 
 times. 
The required solution is then given by
\begin{align}
\underline \psi_{\sfm} =\, {\=M} \,{\underline \sfZ}.
\end{align}
While this is more convoluted than simply performing the tensor multiplication,
this matrix multiplication method is comparatively very fast because  $\~\Psi^{0,0,1}$
is independent of the solution (i.e., constant) and so can be precomputed.

\subsection{Upwinding}\label{upwind}
One of the reasons for choosing {the HDG method} is  to include upwinding in a natural way through the choice of hybridization functions $\tau_m$~\cite{FuQiuHDG}. Upwinding was included in {our simulations} as follows. 
 
{First, the dimensionless effective average electron speed $u_n(z)=-\Jn(z)/n(z)$ and
 the effective average hole speed $u_p(z)=\Jp(z)/p(z)$ were formulated. Next, Eqs.~(\ref{Jn}) and (\ref{Jp}) were used  to obtain 
\begin{align}
u_n(z) & = \mun(z) \frac{d}{dz}\left\{\phi(z)+\phin(z)-\ln \left[n(z)\right]\right\}
\\
\intertext{and}
u_p(z) & = -\mup(z) \frac{d}{dz}\left\{\phi(z)+\phip(z)+\ln \left[p(z)\right]\right\}\,,
\end{align}
respectively. The terms containing $\ln[...]$}  are identifiable as the diffusion terms, with the remaining terms modeling drift in the effective electric field. In drift-dominated regions of the solar cell, which often {constitute most of $\Omega_{\rm el}$,} we therefore get
\begin{align}
u_n(z) & =-\lambda^{-2}\mun(z)  \left[\Edc(z)+E_n(z) \right] 
\\
\intertext{and}
u_n(z) & =\lambda^{-2}\mup(z)  \left[\Edc(z)+E_p(z) \right] \,,
\end{align}
where 
\begin{equation}
E_n(z)=-\lambda^2\frac{d}{dz}\phin(z)
\end{equation}
is the dimensionless effective electric field acting on electrons and
\begin{equation}
E_p(z)=-\lambda^2\frac{d}{dz}\phip(z)
\end{equation}
is the dimensionless effective electric field acting on holes. These two fields arise from  material nonhomogeneity.

 Information in this {drift-dominated} system travels in the {directions} of the electron and hole velocities. At each node, we therefore want to take values of $n(z)$ and $p(z)$ from the inflow side and pass them to the outflow side. As $\mun(z) > 0$ and $\mup(z) > 0$, the respective directions entirely depend on the signs of $\Edc(z) + E_n(z)$ and $\Edc(z) + E_p(z)$. While the hybridization function
 $\tau_m(z)$ is defined as a function of position, only the values at the ends of the elements contribute to the model. Consequently, we {create $\bar{\tau}(z)$} as a piecewise linear function, defined by the limiting values at each side of each node. Two possible limiting values {$\bar{\tau}_1$ and $\bar{\tau}_2$ are chosen such that $\bar{\tau}_1 \ll\bar{\tau}_2$.} Then, to ensure that the correct information flow is achieved,  
\begin{itemize}
\item
$\tau_m\vert_{\gamma,R}=\bar{\tau}_1$ and
$  \tau_m\vert_{\gamma,L} =\bar{\tau}_2$   
are used
if $\Edc(z) +E_n(z)> 0$, or
\item
$\tau_m\vert_{\gamma,L} =\bar{\tau}_1$ and 
$\tau_m\vert_{\gamma,R}=\bar{\tau}_2$  are used
 $\Edc(z) +E_n(z)< 0$.
\end{itemize}
Our choices of $\bar{\tau}_1$ and $\bar{\tau}_2$ are given in Table~\ref{Tab:OptimParam}. 

\subsection{Recombination}\label{numrecom}
The recombination term $R(n,p;z)$ is a nonlinear function of $n(z)$ and $p(z)$. Consequently, care needs to be taken on how to incorporate this term into the HDG method. 
For example, if radiative recombination given by Eq.~(\ref{rad}) is incorporated in  the HDG method, we get
\begin{align}
\pt{R_{rad},\psi_j} 
& = n_\ell p_k \pt{\frac{\alpha_{rad}}{\bar{n}^2} \psi_\ell \psi_k ,\psi_j} - \pt{\alpha_{rad},\psi_j} .
\end{align}
The first term on the right side is nonlinear in the basis function, giving a third-rank tensor with indices $j$, $k$ and $\ell$. Note that $\alpha_{rad}$ and $\bar{n}$ are also projected onto $\mathbb{V}_h$ but, because both are material properties (and are therefore independent of the solution), they do not increase the rank of the tensor produced. The term can thus be implemented in a similar manner to the nonlinear drift term discussed in Section~\ref{Sec:Nonlinear}.

The SRH recombination is given by
Eq.~(\ref{SRH}). If we test this term against the standard polynomial basis $\psi_j$ and integrate over {$\Omega_{\rm el}$, we cannot write the result as a tensor,  because the SRH term is not a polynomial in the basis functions.} Consequently, we integrate using quadrature. For this, we chose Gauss--Lobatto quadrature \footnote{The interpolation points for
the finite-element basis functions are   the Gauss--Lobatto points, as discussed in Sec.~\ref{HDGD}.} to maximize the order of the polynomials that can be exactly integrated, while maintaining use of the values at the nodes. In particular, for a general function $f$, the use of Gauss--Lobatto quadrature to  integrate over {the element} $s_\gamma$ gives
\begin{align}
{\int_{s_\gamma}} 
f[n(z), p(z)] \psi_j(z) \,{\rm d}z &= {\int_{s_\gamma}} 
f\left[\sum_{\ell=1}^{\Ideg} n_\ell \psi_\ell(z), \sum_{\ell=1}^{\Ideg} p_\ell \psi_\ell(z) \right] \psi_j(z) \,{\rm d}z\\
&\approx w_j f\left(n_j, p_j \right),
\end{align}
where $\{w_j\}_{j=1}^{\Pdeg+1}$ are the Gauss--Lobatto quadrature weights for the element $s_\gamma$ and $\Ideg$ is the degree of integration.  A computationally quick way is to use $\Ideg = \Pdeg+1$ quadrature points per element~\cite{masslumping,DDinterp} as this does not require interpolation of the solution.

The Auger recombination is given by Eq.~(\ref{Aug}). As with the SRH term, we cannot write this as a tensor. Either  mass-lumping ~\cite{masslumping} or quadrature \cite{str73} must then be used to calculate the integrals formed when the Auger term is incorporated into the HDG method.

\subsection{Heterojunctions}\label{numhetero}
To implement the jump conditions (\ref{hn}) or (\ref{hp}), we ensure that a node falls at each discontinuity in  the material parameters. Every jump in $n(z)$ can be taken into account by redefining the jump operator as
\begin{align}
\dsbg{\tau} \hat{n}_\gamma &= \left[\tau\vert_{\gamma+1,L} \exp(\Delta \chi_{\gamma}) - \tau\vert_{\gamma,R}\right]  \hat{n}_\gamma
\intertext{and the difference operator as}
\ipt{\hat{n}, \zeta \psi_j} &= (\zeta\psi_j)\vert_{\gamma,R} \hat{n}_{\gamma} - (\zeta\psi_j)\vert_{\gamma,L} \exp(\Delta \chi_{\gamma-1})\hat{n}_{\gamma-1}\,,
\end{align}
where $\Delta \chi_{\gamma}$ is the jump in electron affinity $\chi$ across node $\gamma$ and $\zeta$ is a place-holder function. This is equivalent to defining to values of $\hat{n}_\gamma$ at each node, separated in value by the necessary jump induced by the discontinuous electron affinity. Jumps in $p(z)$ can be handled similarly.

\subsection{Homotopy}\label{homotopy}
In order to aid convergence of the highly nonlinear discrete system, homotopy  is employed in our
simulation. The fixed charge density $\Nf$, recombination rate $R(n,p;z)$, and bandgap non-homogeneity, i.e. $E_g(z) - E_{g,av}$ where $E_{g,av}$ is the mean bandgap, are all multiplied by a constant $\delta_{\rm homo}$. The simulation is started with $\delta_{\rm homo}^{(0)} = \delta_{\rm min}$. Once    the simulation is deemed to have converged,   a larger value 
 \begin{align}
\delta_{\rm homo}^{(1)} = \text{min}(\delta_{\rm homo}^{(0)} \eps_{\rm homo}, 1)\,
\end{align}
is chosen for $\delta_{\rm homo}$,
with {$\eps_{\rm homo}>1$} as the homotopy damping relaxation rate. 
Thus, the iteration
\begin{align}
\delta_{\rm homo}^{(\ell+1)} = \text{min}(\delta_{\rm homo}^{(\ell)} \eps_{\rm homo}, 1)\,
\end{align}
is employed, full convergence being deemed to occur
in the iteration in which     $\delta_{\rm homo}$ becomes equal to unity.

\section{Numerical Test of HDG Method}\label{HDGtest}
As a model test problem to show the behavior of the HDG {method (but not to present a solar-cell design), the} solar cell was taken to comprise a p-i-n junction made from copper-indium-gallium-(di)selenide (CIGS), {as shown in Fig.~\ref{Ben1}(a).} 
The electrical parameters of CIGS  \cite{Frisk14} are given in Table~\ref{Tab:CIGS}. For this test
problem, the generation rate was taken to be uniform, with $\JscOpt = 10$~mA~cm$^{-2}$; thus RCWA results were not used.  The $p$-type layer of thickness $\Lp=20$~nm was taken to be doped with acceptor atoms with concentration specified via $\Nf = -10^{17}$~cm$^{-3}$, and the 
$n$-type layer of thickness $\Ln=20$~nm  with donor atoms with concentration specified via $\Nf = 10^{17}$~cm$^{-3}$. The undoped $i$-type layer was taken to be $\Li=200$~nm thick;
{thus, $\Lz=\Lp+\Li+\Ln=240$~nm. All three layers were taken to be homogeneous. Neither the radiative nor the Auger recombination mechanism was activated for the test problem.}

\begin{table}
	\centering
\begin{tabular}{|c|c|c|}\hline
Parameter & Symbol & Value\\\hline\hline
Bandgap  & $E_g$ & 1.3~V\\
Electron affinity & $\chi$ & 4.5~V\\ 
Conduction-band density of states & $\Nc$ & $2.22\times 10^{18}$~cm$^{-3}$\\
Valence-band density of states & $\Nv$ & $1.78\times 10^{19}$~cm$^{-3}$\\
DC relative permittivity & $\eps_{dc}^0$ & 13.6\\
Electron mobility & $\mun$ & $100$~cm$^2$~V$^{-1}$s$^{-1}$\\
Hole mobility & $\mup$ & 25~cm$^2$~V$^{-1}$s$^{-1}$\\
SRH lifetime parameter (electrons)&$\taun(z)$&$1\times10^{-9}$s$^{-1}$\\
SRH lifetime parameter (holes)&$\taup(z)$&$1\times10^{-9}$s$^{-1}$\\\hline
\end{tabular}
\caption{Electrical parameters of CIGS.}\label{Tab:CIGS}
\end{table}

Chiefly, three    parameters  affect the accuracy of the HDG method: the degree of the interpolating polynomials {$\Pdeg$}, the length of each element $d_z$, and the degree of quadrature integration $\Ideg$. The convergence of the {method with respect to each of these parameters was} investigated. The parameters  used for the convergence study are given in Table~\ref{Tab:HDGsetup}.

\begin{table}
\centering
\begin{tabular}{|c|c|c|}
	\hline
Symbol 		& Value 	& Description\\\hline\hline
$\Pdeg$ 	& variable 		& Degree of {Lobatto polynomials}\\
$N_{\rm zp}$		&  {variable} 		& Number of points in $p$-type layer\\
$N_{\rm zi}$		&  variable 		& Number of points in $i$-type layer\\
$N_{\rm zn}$		&  variable 		& Number of points in $n$-type layer\\
$\Ideg$		& variable		& Nonlinear term integration degree\\\hline
${\bar{\tau}}_1$			& $10^{-3}$ & First parameter for hybridization\\
${\bar{\tau}}_2$			& $10^{3}$ 	& Second parameter for hybridization  \\
$\delta V$	& 0.01		& Maximum  step for {$\Vext$} \\
$\delta_{\rm min}$		& $10^{-2}$			& Initial homotopy damping\\
{$\eps_{\rm homo}$} 	& $1.2$		& Homotopy damping relaxation\\
$\delta_{\rm num}$		& $1$			& No. of homotopy attempts\\
$\delta_{\rm inc}$	& 10		& Homotopy increase on new attempt\\
$n_{\rm loop}$	& 10		& {Maximum number of  iterations for Newton's method}\\
$E_{\rm rel,tol}$		& $10^{-4}$	& Relative change in state vector allowed\\
$E_{\rm abs,tol}$		& $10^{-6}$	& Absolute change in state vector allowed\\
$J_{\rm tol}$ 	& $0.1$		& Relative noise allowed in $J$\\
$P_{\rm ref}$ 	& 100		& Max iteration to find $\Pmax$\\
$P_{\rm tol}$		& $10^{-6}$	& Allowed variation around $\Pmax$\\\hline
\end{tabular}	
\caption{Baseline parameter settings in Sec.~\ref{HDGtest} for the HDG simulation of the model CIGS p-i-n solar cell  in Fig.~\ref{Ben1}(a). {Note
that $\Nz=N_{\rm zp}+N_{\rm zi}+N_{\rm zn}$. {The \textit{state vector} is the vector of all solution values in the electronic step.}}
}\label{Tab:HDGsetup}
\end{table}

The convergences of
the short-circuit current density $\Jsc$, the open-circuit voltage $\Voc$, the maximum power  $\Pmax$, and the fill factor $FF$  with respect to $\Pdeg$, $\dz$, and $\Ideg$, were investigated.
 Figure~\ref{Fig:HDGconva}  shows the {error in each of four electrical characteristics relative to the value
 for $\Pdeg=9$ when  $\dz= 2$~nm and $\Ideg=10$.}
As  $\Pdeg$ increases from $2$, each of the four errors is seen to decrease exponentially, with a rate approximately proportional to $\exp(-3\Pdeg)$. This reduction in error is seen to saturate at around $\Pdeg = 5$.
 
\begin{figure}
\begin{center}
\includegraphics[width=0.6\textwidth]{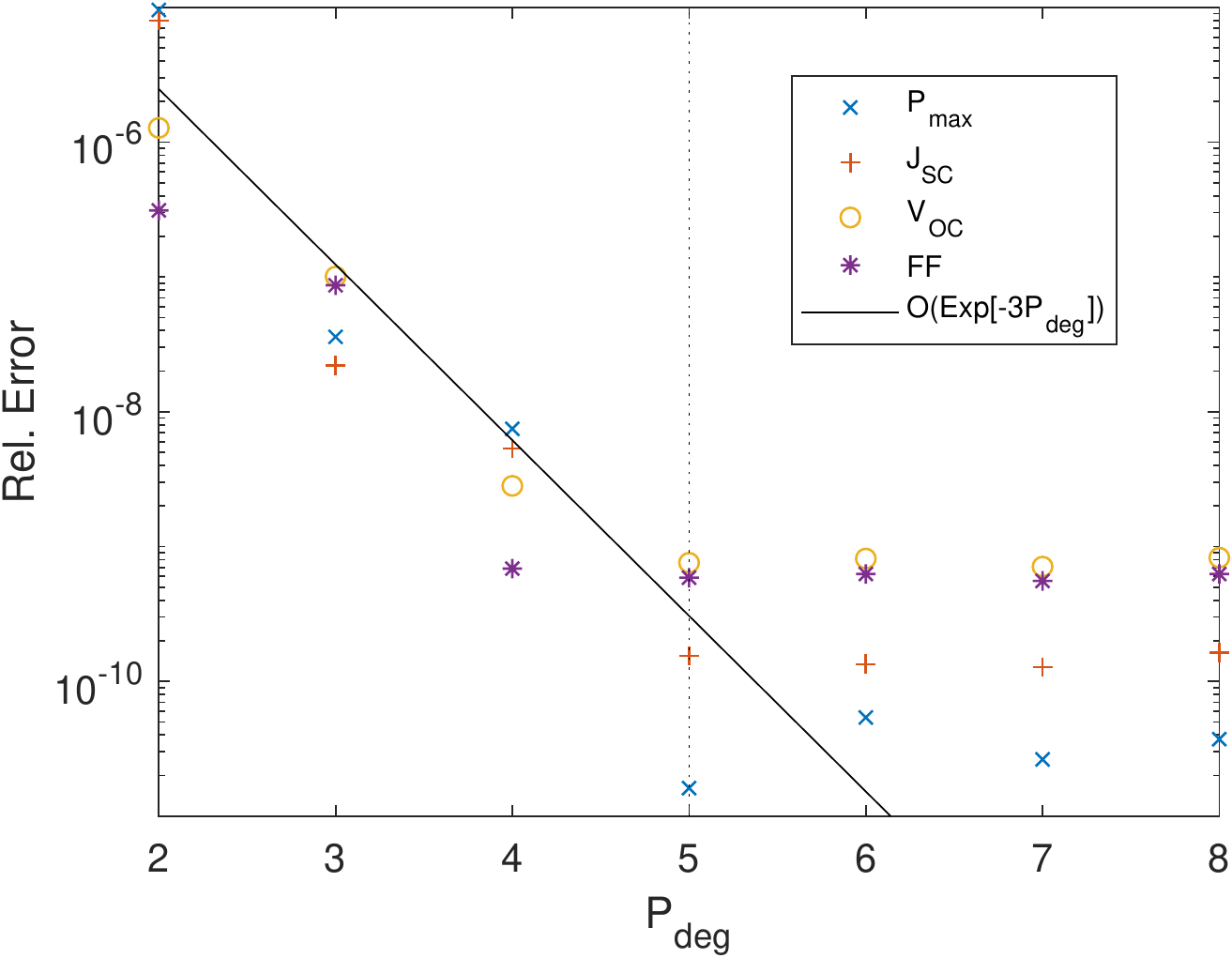}
\end{center}
\caption{Relative
errors in $\Pmax$, $\Jsc$, $\Voc$,    and  $FF$ against  $\Pdeg\in[2,8]$,
when {$\dz= 2$~nm and $\Ideg = 10$.}
}\label{Fig:HDGconva}
\end{figure}

Figure~\ref{Fig:HDGconvb}   shows the {error in each of four electrical characteristics relative to the value for $\dz=1$~nm when  {$\Pdeg=5$}  and $\Ideg=10$. 
As  $\dz$ decreases from $20$~nm, all four errors  decrease as $\mathcal{O}(\dz^4)$ {which is suboptimal}.}

\begin{figure}
\begin{center}
\includegraphics[width=0.6\textwidth]{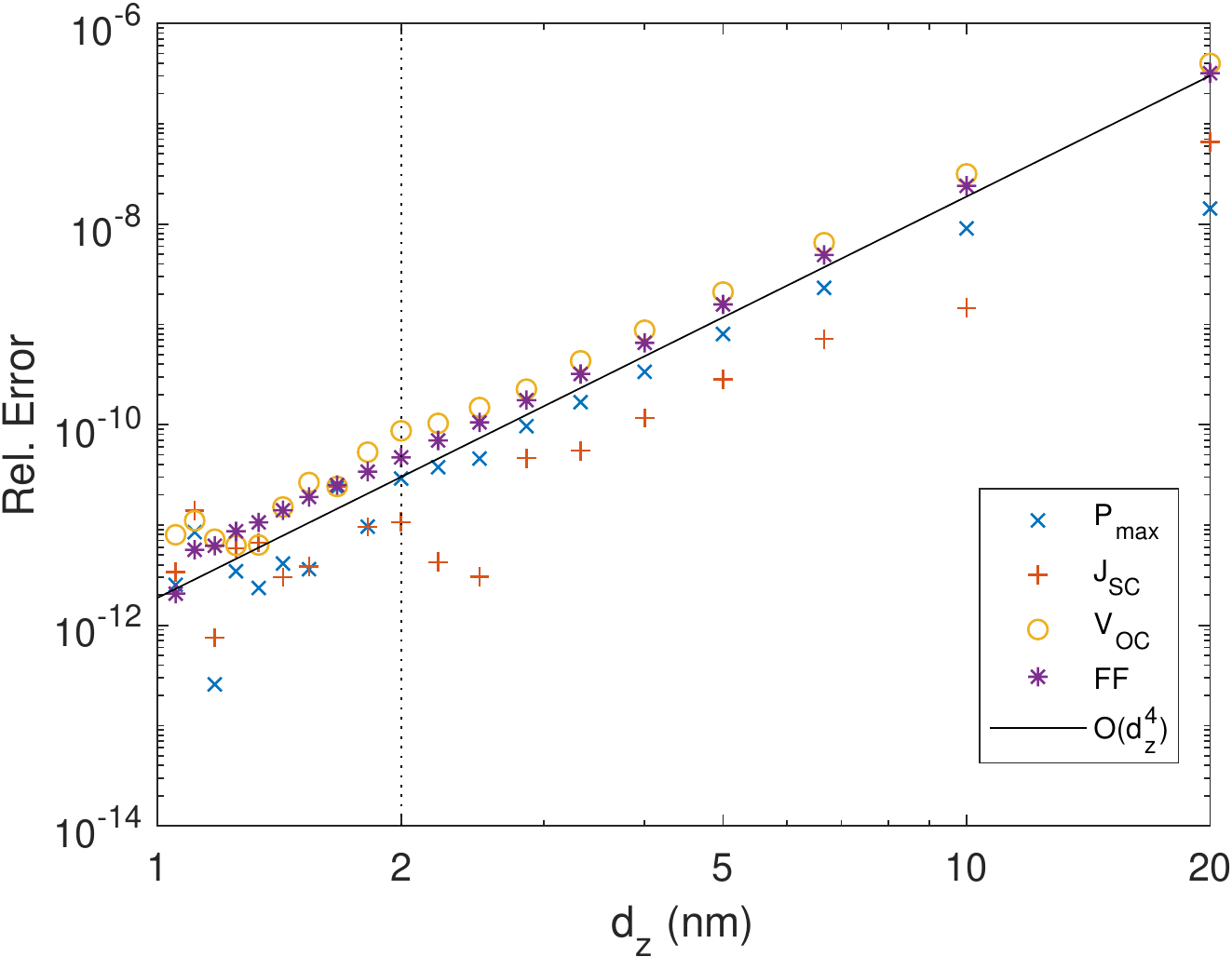}
\end{center}
\caption{Relative
errors in $\Pmax$, $\Jsc$, $\Voc$,    and  $FF$ against {$\dz\in[1,20]$~nm,
when  $\Pdeg= 5$ and $\Ideg=10$.}
}\label{Fig:HDGconvb}
\end{figure}

Finally, Fig.~\ref{Fig:HDGconvc} shows the {error in each of four electrical characteristics relative to the value
 for $\Ideg=12$ when  $\Pdeg=5$ and $\dz=2$~nm.}
As  $\Ideg$ increases, each error is seen to decrease exponentially, with a rate approximately proportional to $\exp(-4\Ideg)$. This reduction in error is seen to saturate at around $\Ideg = 5$.

\begin{figure}
\begin{center}
\includegraphics[width=0.6\textwidth]{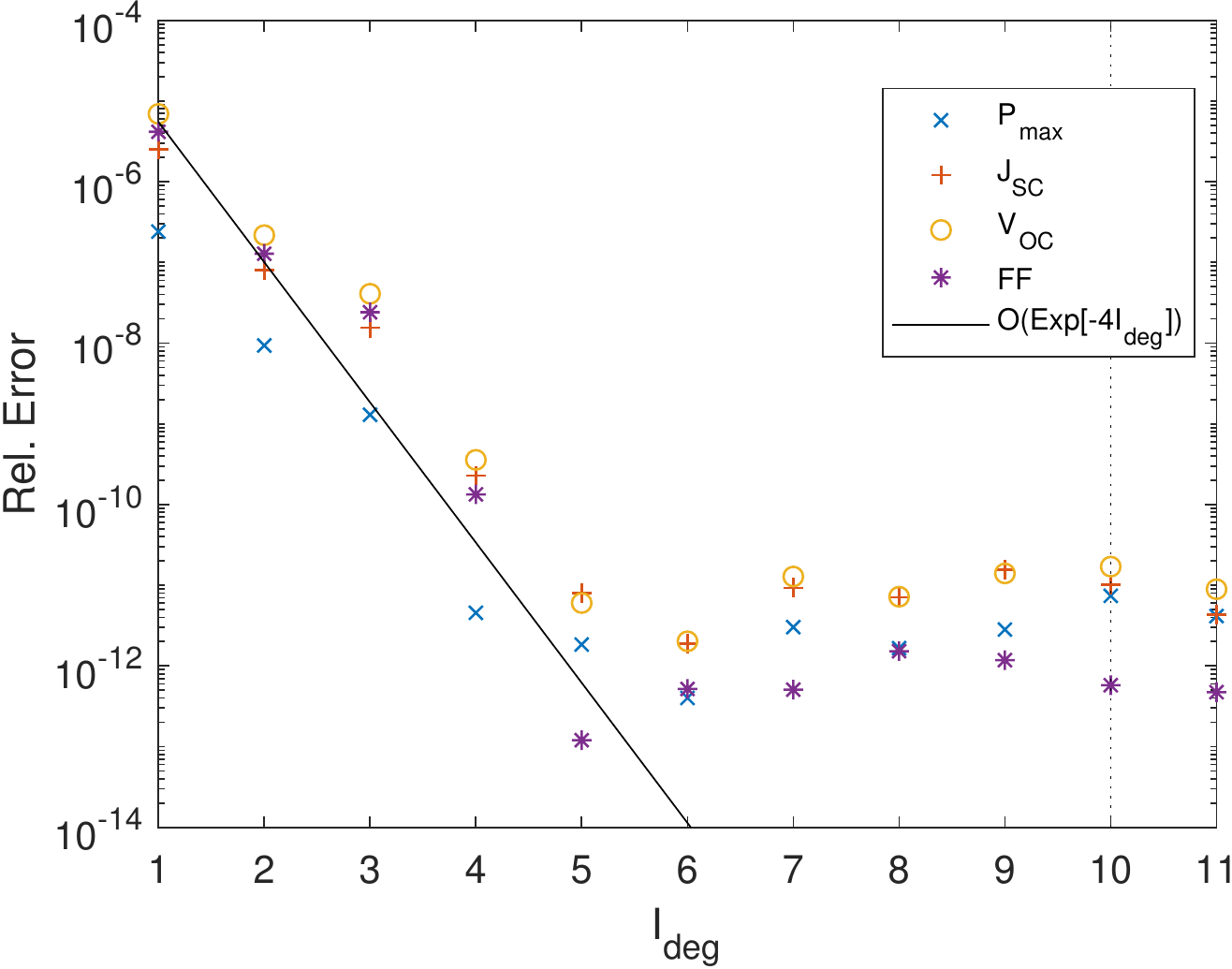}
\end{center}
\caption{{Relative
errors in $\Pmax$,  $\Jsc$, $\Voc$,   and  $FF$ against  $\Ideg\in[1,11]$,
when {$\Pdeg=5$ and $\dz= 2$~nm.}}
}\label{Fig:HDGconvc}
\end{figure}

Thus,
Figs.~\ref{Fig:HDGconva}--\ref{Fig:HDGconvc} {show that  the short-circuit current density $\Jsc$, the open-circuit voltage $\Voc$, the maximum power  $\Pmax$, and the fill factor $FF$, all converge at approximately the same rate with respect to $\Pdeg$, $\dz$, and $\Ideg$}.
The order of convergence is not optimal, as might be expected {since we used polynomials of the same degree in all elements.}

\section{Differential Evolution Algorithm for Optimization}\label{DEA}
The differential evolution algorithm \cite{DEA} is a gradient-free optimization algorithm suited for maximizing an objective function over a high-dimensional parameter space. Given the objective function $C: {\mathbb{S}}\to\mathbb{R}$, where $ {\mathbb{S}}\subset \mathbb{R}^{\tilde{N}}$ is the set of all possible choices of the $\tilde{N}$ input parameters, the DEA starts by selecting $\NP$ points to form an initial population $\#P_0 \subset {\mathbb{S}}$. The objective function $C$ is evaluated at each of these points, with the results used by DEA in a  mutation-recombination-selection process to build a new population $\#P_1 \subset {\mathbb{S}}$. The objective  function is then evaluated at {each point of} this new population to develop a new population $\#P_2 \subset {\mathbb{S}}$ in the second  mutation-recombination-selection step.  This process continues
iteratively until the objective function appears to stabilize.

The  efficiency $\eta$ defined in Eq.~(\ref{etadef}) is the appropriate objective
function for designing solar cells. Following suggestions 
in the DEA documentation~\cite{DEA}, we fix  the crossover fraction as $C_{\rm R} = 0.6$. The step size $F$ used in the  mutation steps is set to be randomly distributed in $\left[0.5,1\right]$ uniformly. Allowing $F$ to vary randomly with each iteration has been termed \textit{dither}, and its use has been shown to improve convergence for many {problems \cite{Swagatam2011}.}

\section{Numerical Test of DEA}\label{DEA-test}

As a model test problem,  solar cell was taken to comprise a p-i-n junction made from CIGS, as shown in Fig.~\ref{Ben1}(a).
The electrical parameters of CIGS  \cite{Frisk14} are given in Table~\ref{Tab:CIGS}.   
The $p$-type layer of thickness {$\Lp=20$~nm was taken to be doped with acceptor atoms with concentration specified via $\Nf = -10^{17}$~cm$^{-3}$, and the 
$n$-type layer of thickness $\Ln=20$~nm  with donor atoms with concentration specified via $\Nf = 10^{17}$~cm$^{-3}$.} The thickness $\Li$ of the undoped $i$-type layer was taken as a variable
for optimization. All three layers were taken to be homogeneous, with bandgaps $\sfE_{\rm g,p}$, $\sfE_{\rm g,n}$, and
$\sfE_{\rm g,i}$ taken to be variables for optimization.
Neither the radiative nor the Auger recombination mechanism was activated for the test problem.

{The antireflection coating in Fig.~\ref{Ben1}(a) was ignored. We fixed $\Lm=150$~nm. The period
$\Lx$ was kept variable. The height $\Lg$ and base $\zetag\Lx$, $\zetag\in(0,1)$, of the metallic protuberance  in each period were also kept variable.}
The $\tilde{N}=7$ parameters   chosen to be optimized in our test problem thus are: $\Lx$, $\Lg$, $\zetag$,
$\Li$, $\sfE_{\rm g,p}$, $\sfE_{\rm g,n}$, and
$\sfE_{\rm g,i}$. 

Both the photonic and the electronic steps were implemented. The population number was set as $\NP = 100$. 
The remaining parameters used for the optimization algorithm were as described in the previous sections, together with values given in Table~\ref{Tab:OptimParam}. These values were chosen to maintain reasonable accuracy, aid convergence for a wide range of solar cell designs, and to allow rapid computation of {the efficiency at each step of DEA.} {These choices} resulted in very quick    evaluation of $\eta$ at a particular choice of the parameter values in roughly 6~min, with four evaluations running concurrently in
MATLAB on a 20-processor (Intel Xeon Gold 6138\@2GHz) Linux (Ubuntu 17.10) computer.  The actual time per parameter set depends on the HDG solver (being longer if more homotopy steps
are needed).

\begin{table}
\centering
\begin{tabular}{|c|c|c|c|}
	\hline
 Symbol 		& Value 	& Description\\\hline
 $\Nt$ 	& 10 		& Fourier order for RCWA \\
 $\Pdeg$ 	& 3 		& {Degree of Lobatto polynomials}\\
 $N_{\rm zp}$		&  {10} 		& Number of points in $p$-type layer\\
$N_{\rm zi}$		&  {50} 		& Number of points in $i$-type layer\\
$N_{\rm zn}$		&  {10} 		& Number of points in $n$-type layer\\
 $\Ideg$		& 4		& Nonlinear term integration degree\\\hline
  $\bar{\tau}_1$			& $10^{-6}$ & First parameter for hybridization  \\
 $\bar{\tau}_2$			& $10^{6}$ 	& Second parameter for hybridization  \\
$\delta_{\rm min}$		& $10^{-2}$			& Initial homotopy damping\\
{$\eps_{\rm homo}$} 	& 1.1		& Homotopy damping relaxation\\
$\delta_{\rm num}$		& 40			& No. of homotopy attempts\\
 $\delta_{\rm inc}$	& 5		& Homotopy increase on new attempt\\ 
$n_{\rm loop}$	& 40		& {Maximum number of iterations for Newton's method}\\
{$E_{\rm rel,tol}$}		& $10^{-2}$	& Relative change in state vector allowed\\
{$E_{\rm abs,tol}$}		& $10^{-4}$	& Absolute change in state vector allowed\\
 $P_{\rm tol}$		& $10^{-3}$	& Allowed variation around $\Pmax$\\ 
  $t_{\rm max}$ & 300~s & Maximum computer time \\
  \hline
\end{tabular}	
\caption{Baseline parameter settings for DEA optimization in Sec.~\ref{DEA-test}. These are laxer than those used in Sec.~\ref{HDGtest} as rapid computation is necessary. Other settings are the same as in Table~\ref{Tab:HDGsetup}.}
\label{Tab:OptimParam}
\end{table}

Figure~\ref{FIG:DEAConv} summarizes the progression of the DEA towards an
optimal result.  Each of panels (a)-(g) in this figure is
a graph of the efficiency $\eta$ as a function of one of the parameters in the optimization. Thus,
each point corresponds to a 7-dimensional vector of parameter values used by DEA.  The optimal efficiency found is marked by a large red disk.  The remaining panel, Fig.~\ref{FIG:DEAConv}(h), shows the progress of optimization. In particular it shows the efficiency as a function of DEA step.
The  optimal values of the seven parameters in this study are given in Table~\ref{tab_param} which result in  $\eta=15.7\%$. 
\begin{figure}[ht]
\centering
	\begin{tabular}{cc}
		(a)\includegraphics[width=0.45\textwidth]{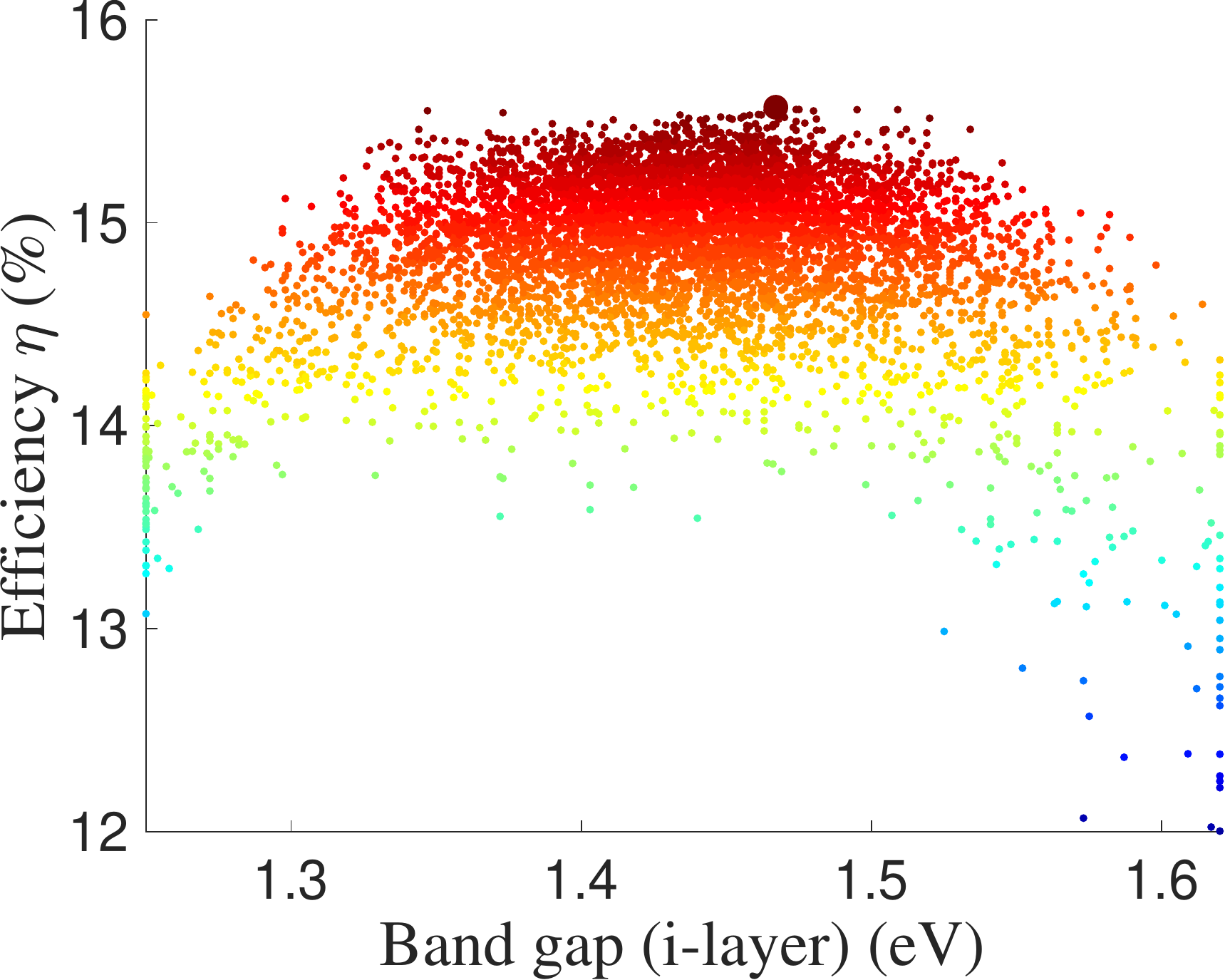}&
		(b)\includegraphics[width=0.45\textwidth]{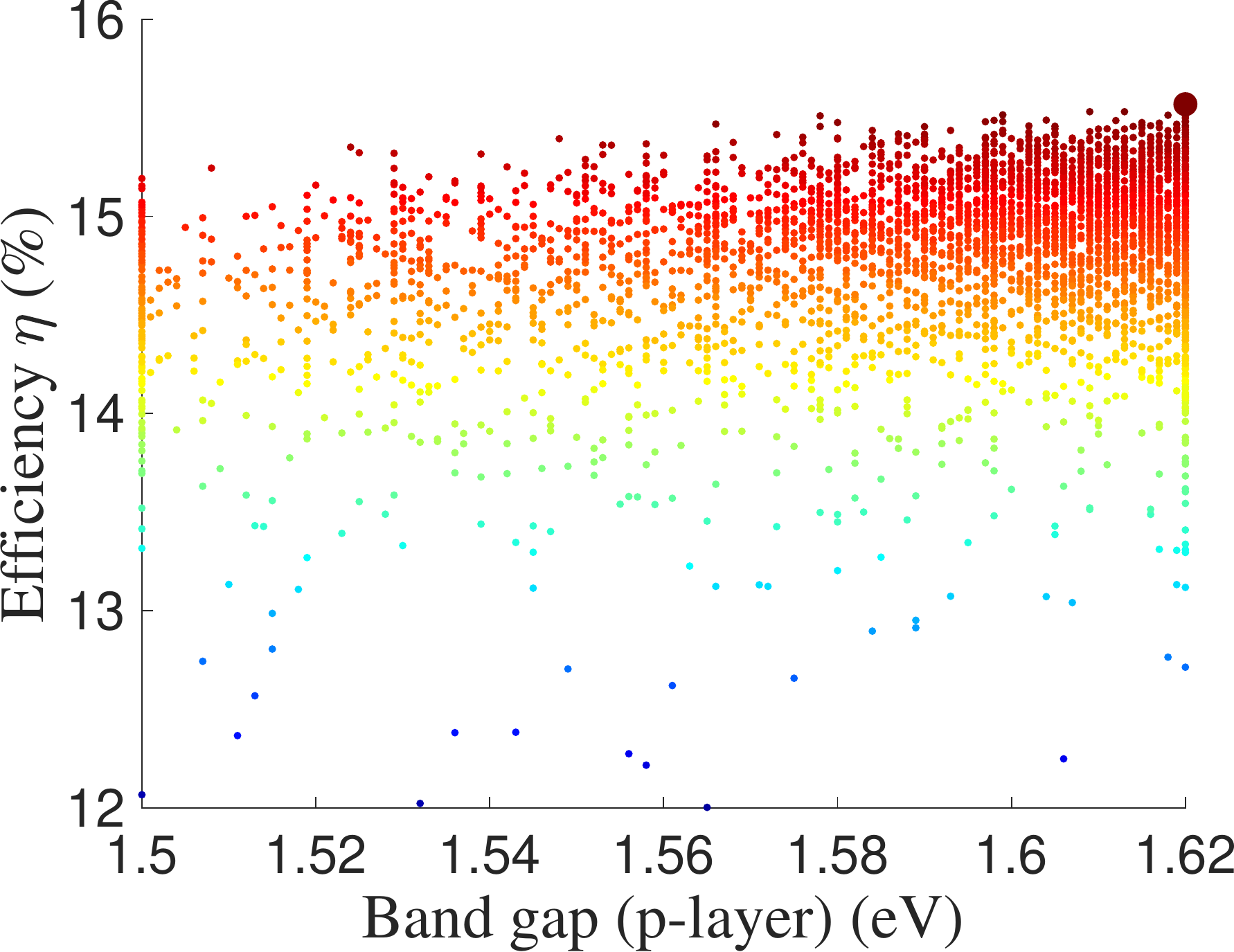}\\
		(c)\includegraphics[width=0.45\textwidth]{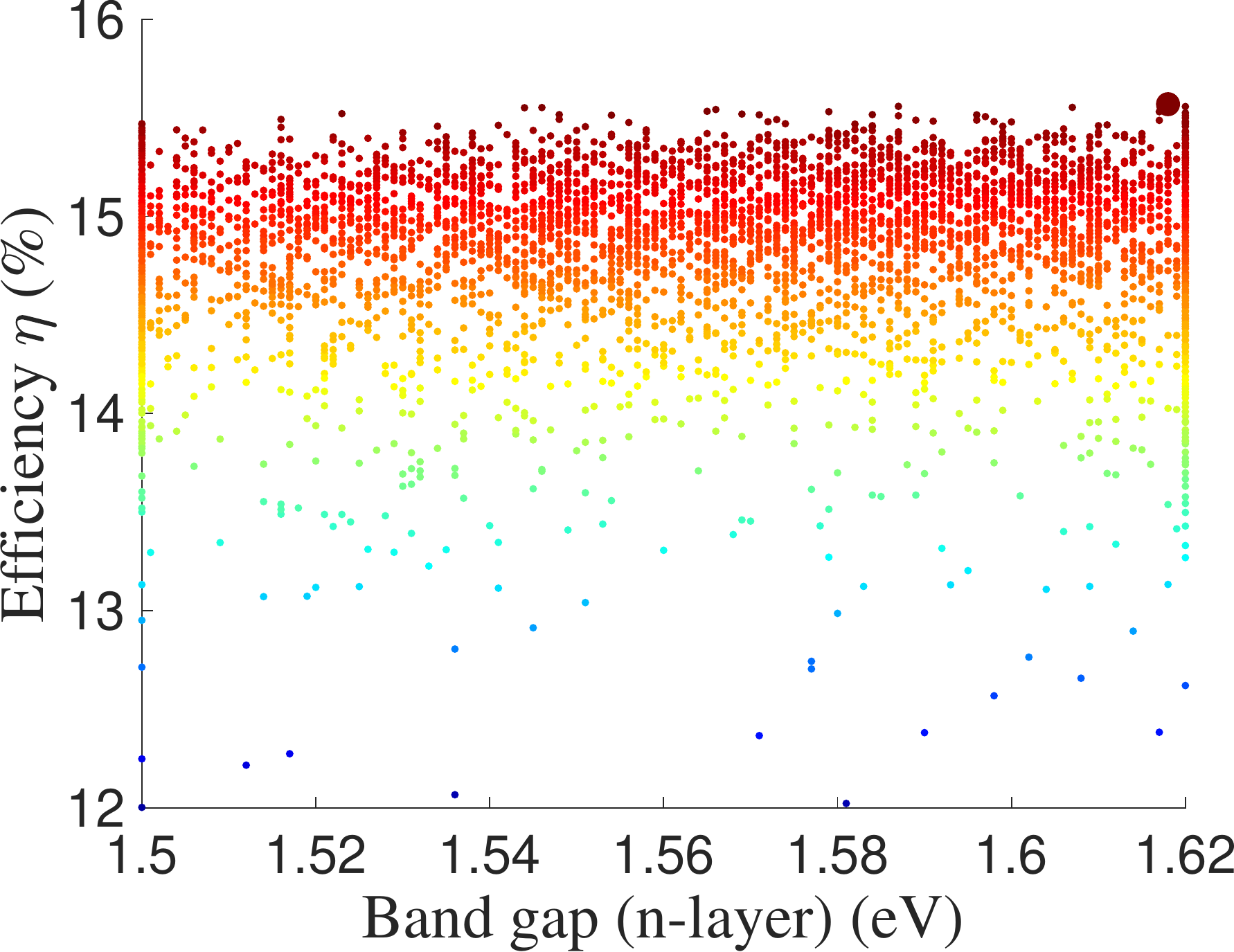}&
		(d)\includegraphics[width=0.45\textwidth]{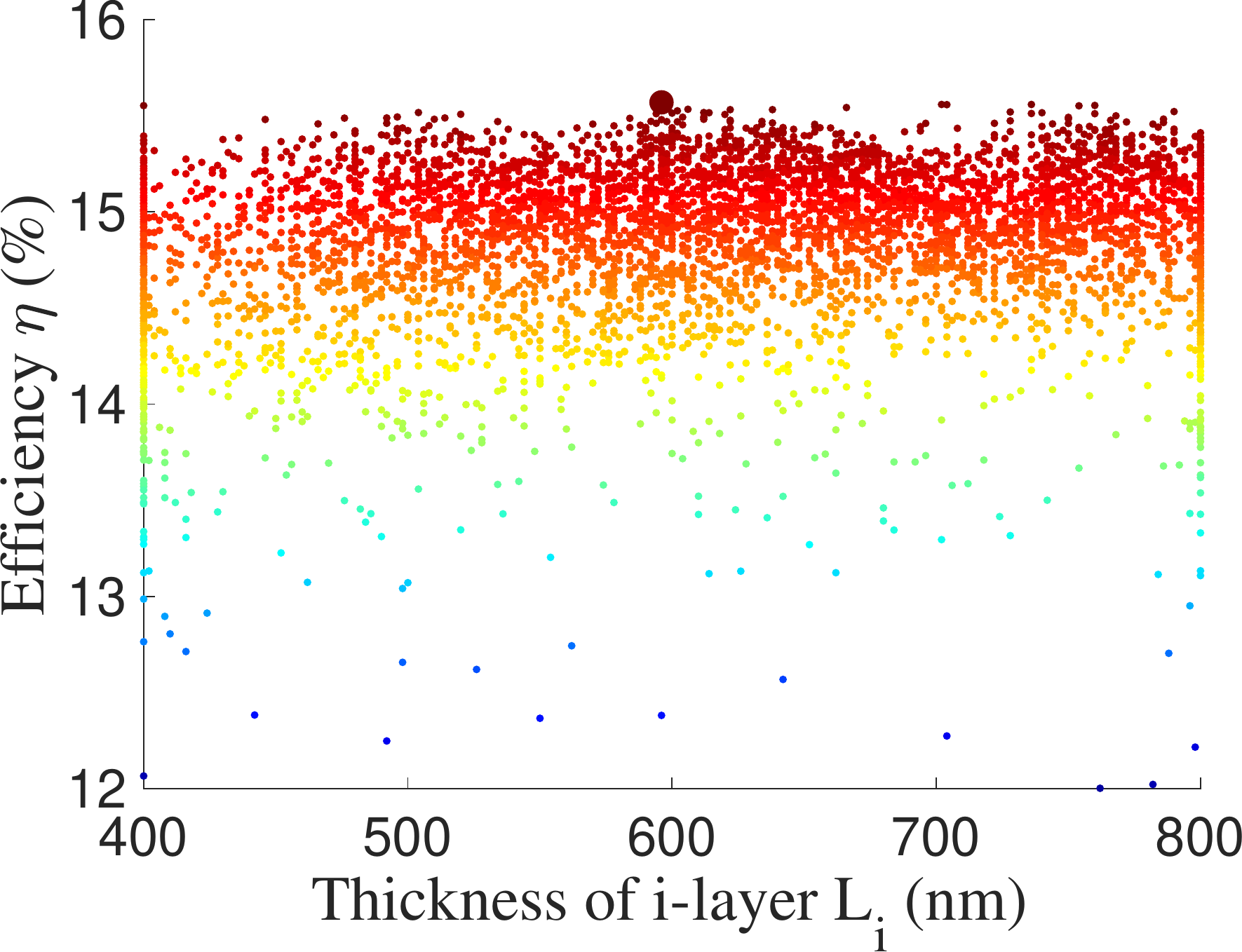}\\
		(e)\includegraphics[width=0.45\textwidth]{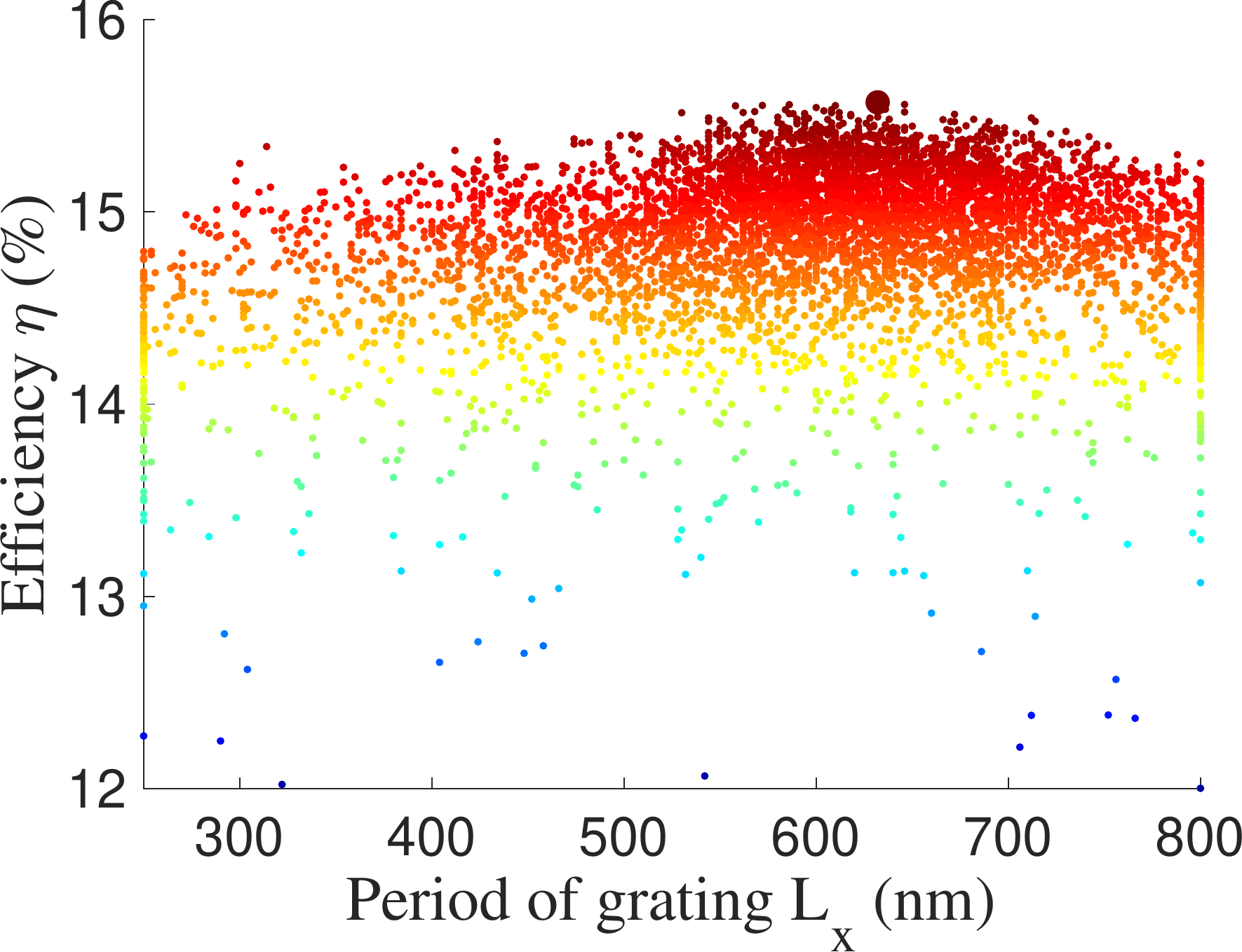}&
		(f)\includegraphics[width=0.45\textwidth]{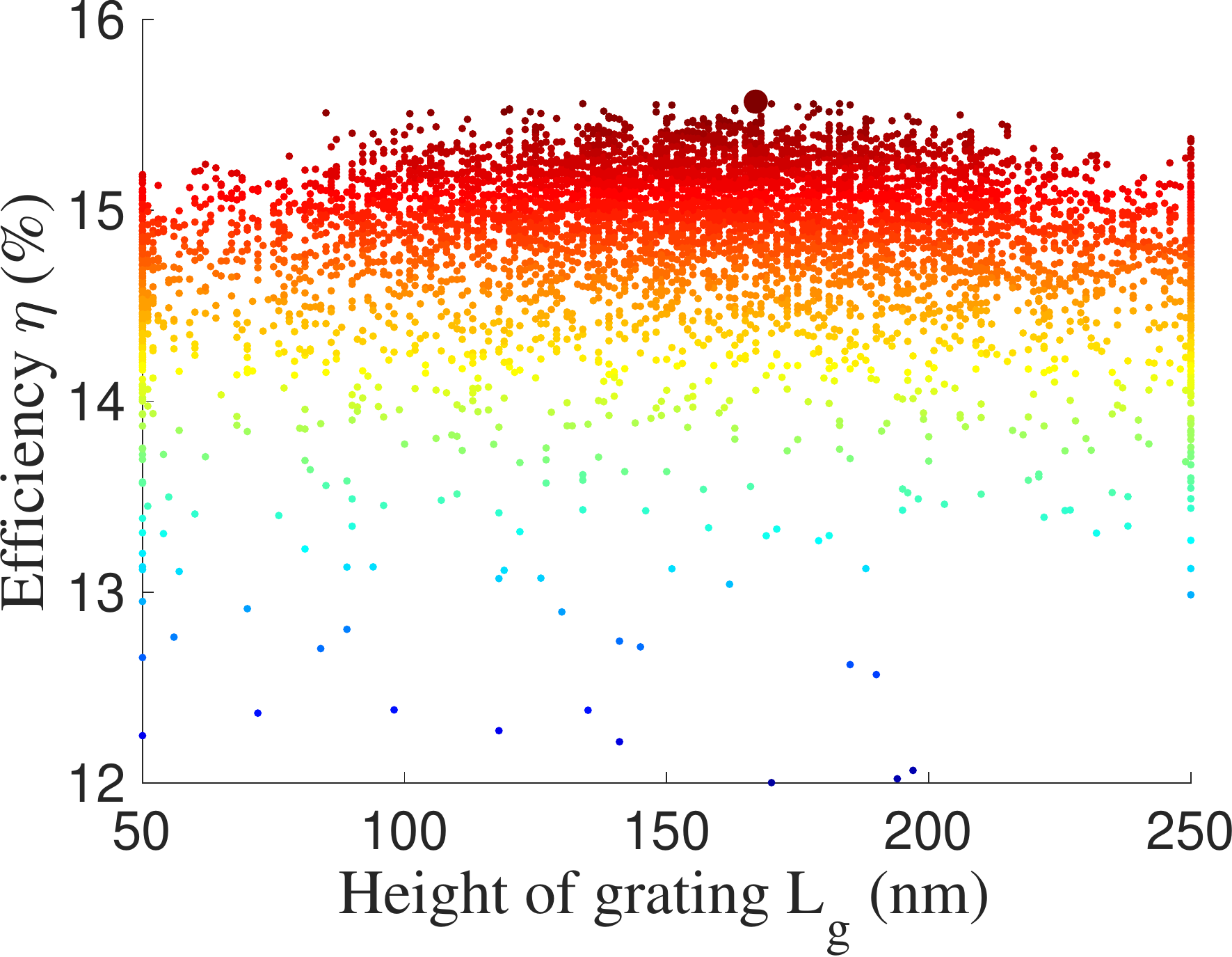}\\
		(g)\includegraphics[width=0.45\textwidth]{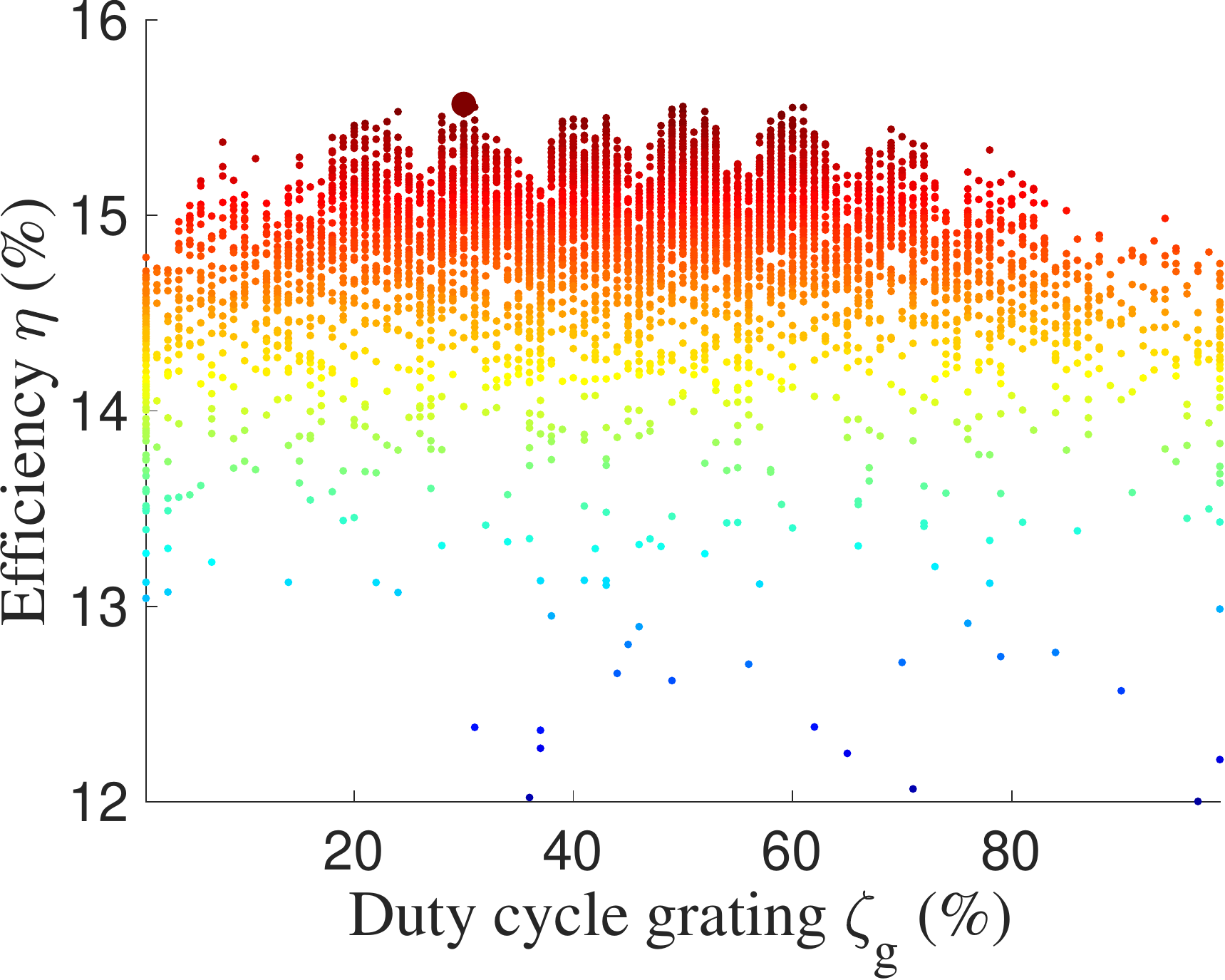}&
		(h)\includegraphics[width=0.45\textwidth]{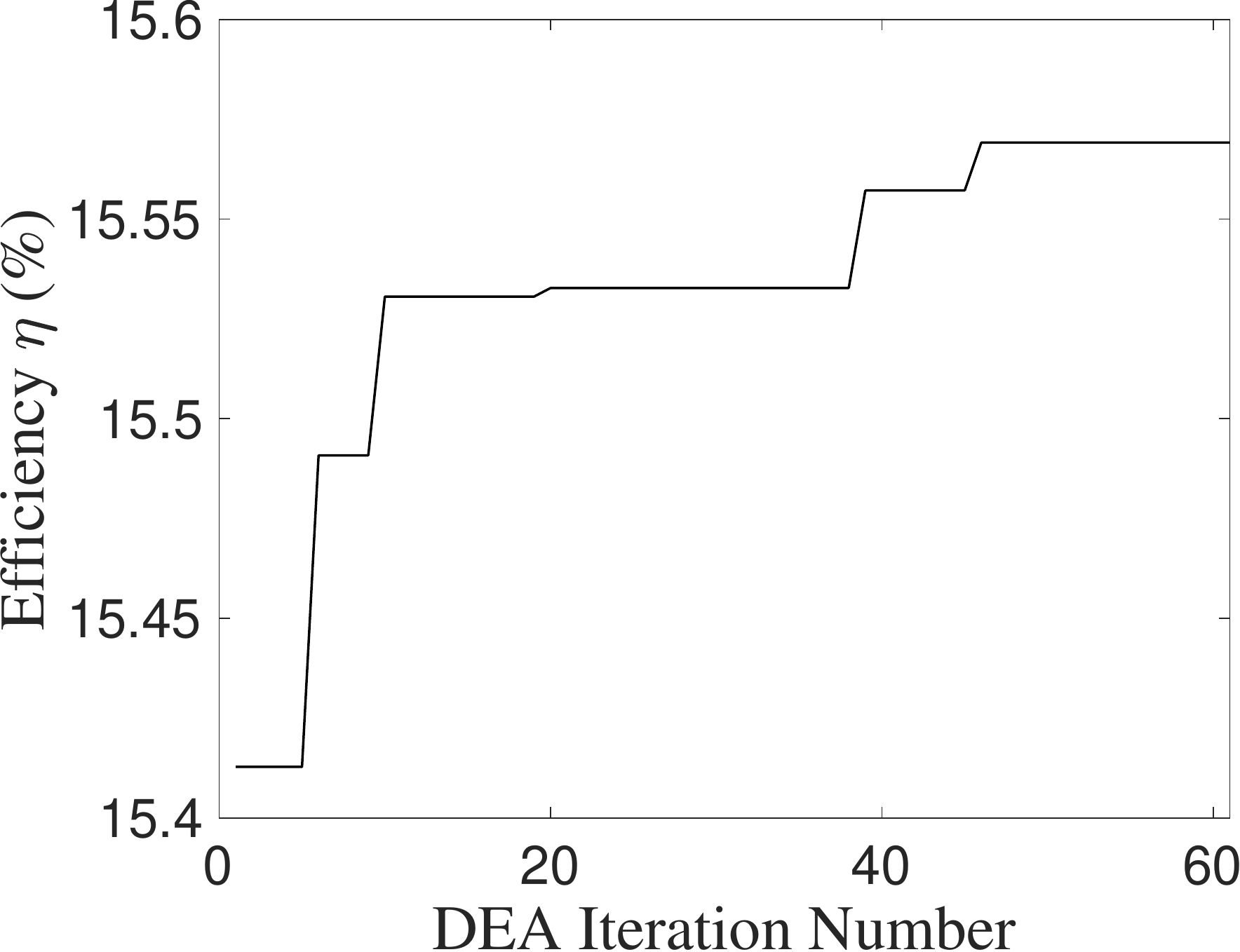}
    	\end{tabular}
\caption{Projections of the results from DEA optimization study of Sec.~\ref{DEA-test}
onto the plane containing $\eta$ and
(a) $\Lx$, (b) $\Lg$, (c) $\zetag$, (d) $\Li$,  (e) $E_{g,p}$,
(f) $E_{g,i}$, (g) $E_{g,n}$,   and (h) DEA iteration number.}\label{FIG:DEAConv}
\end{figure}

\begin{table}[ht]
 \begin{center}
\begin{tabular}{|c|c|}
 \hline
 Variable Parameter &  Optimal Value\\
 \hline\hline
 $\Lx$ &632~nm\\
 $\Lg$ &167~nm\\
 $\zetag$  &{0.3}\\
 $\Li$  &596~nm\\
 $\sfE_{\rm g,p}$  &  1.62~V\\
$\sfE_{\rm g,i}$ &   1.467 ~V\\
$\sfE_{\rm g,n}$ &      1.618~V\\
\hline
 \end{tabular}
 \end{center}
 \caption{Parameters varied in the optimal design study described in Sec.~\ref{DEA}.
 The optimal value refers to the value at the end of the DEA computation.}
 \label{tab_param}
\end{table}

\section{Concluding Remarks}\label{concl}
{We have presented details of a design tool for optimizing 
the efficiency of thin-film photovoltaic solar cells. The solar cell
can have multiple semiconductor layers in addition to dielectric
layers serving as antireflection coatings, passivation layers, and buffer layers.
The solar cell is backed by a metallic grating which is periodic along a fixed
direction. }

{The heart of the design tool is a coupled optoelectronic simulation. The photonic step
of the simulation delivers
the 2D variation of the electron-hole-pair generation rate inside the semiconductor layers
of the solar cell. After averaging along the direction of the periodicity of the grating,
the electron-hole-pair generation rate becomes the input to the electronic step
of the simulation. The chief output of the  electronic step is the efficiency of the solar cell.
The design tool uses the differential evolution algorithm to determine the dimensions
and bandgaps of the semiconductor layers to maximize the efficiency of the solar cell.}

{The design tool can be augmented to incorporate a biperiodic metallic grating at the
cost of increased computation time \cite{Civiletti}. Whether biperiodicity should be incorporated
 will depend on the consequent augmentation of the efficiency of the solar cell, and is therefore
 a matter of further research. Given that the period of the grating is in the 400-to-1000~nm range
 in order to invoke guided-wave phenomena for increased photon trapping \cite{Civiletti,Faiz18,AndersonSPIE,SBS1983,HM1995,Khaleque} but the electronic step is
electrostatic in character,  the averaging of
the electron-hole-pair generation rate delivered by the photonic step  
along the periodicity directions of the grating is appropriate; elimination of that averaging
will significantly increase computational time without significant gain in the electronic step.}

{We have begun to apply the design tool for diverse types of thin-film photovoltaic solar cells
and will report our results in due course of time.}

\section*{Acknowledgements}
A. Lakhtakia thanks the Charles Godfrey Binder Endowment at the Pennsylvania
State University for ongoing support of his research. 
The research of  T.~H. Anderson, B.~J. Civiletti, and P.~B.  Monk was partially supported by  the US
National Science Foundation under grant number DMS-1619904.
The research of A. Lakhtakia was partially supported by  US National Science Foundation  
under grant number DMS-1619901.

\end{document}